\documentclass[12pt]{article}
\usepackage{amsmath,amsthm,amsfonts}

\usepackage{subfigure}
\usepackage[mathscr]{eucal}
\usepackage{graphics,epsfig}
\usepackage{longtable}
\usepackage{algorithmic}
\usepackage[ruled]{algorithm}
\usepackage{multirow}
\usepackage{rotating}
\usepackage{verbatim}

\usepackage{color}

\newtheorem{lem}{Lemma}[section]
\newtheorem{thm}{Theorem}[section]

\newtheorem{cor}{Corollary}[section]
\numberwithin{equation}{section} \numberwithin{figure}{section}
\numberwithin{table}{section}

\newtheorem{rem}{Remark}[section]

\newcommand{\ds}{\displaystyle}

\def\l{\label}    

\def\m{\mbox}   
  \def\ss{\smallskip}

\def\bb{\begin{equation}} \def\ee{\end{equation}}

\def\beqn{\begin{eqnarray}}  \def\eqn{\end{eqnarray}}

\def\beqnx{\begin{eqnarray*}} \def\eqnx{\end{eqnarray*}}

\def\bn{\begin{enumerate}} \def\en{\end{enumerate}}

\def\bd{\begin{description}} \def\ed{\end{description}}

\begin{document}

\begin{center}
{\bf An Inexact Uzawa Algorithm \\
for Generalized Saddle-Point Problems and Its Convergence}
\end{center}

\bigskip
\centerline{Kazufumi ~Ito\footnote{Department of Mathematics, North
Carolina State University, Raleigh, North Carolina, USA. ({\tt
kito@math.ncsu.edu}).}
\quad \quad Hua ~Xiang\footnote{School of Mathematics and
Statistics, Wuhan University, Wuhan, China. ({\tt
hxiang@whu.edu.cn}).}
 \quad \quad Jun ~Zou\footnote{
Department of Mathematics, The Chinese University of Hong Kong,
Shatin, N.T., Hong Kong. The work of this author was substantially
supported by Hong Kong RGC grant  (Project~405110 and 404611). ({\tt
zou@math.cuhk.edu.hk}).} }

\bigskip
\begin{abstract}
We propose an inexact Uzawa algorithm with
two variable relaxation parameters for solving
the generalized saddle-point system. The saddle-point problems
can be found in a wide class of applications, such as
the augmented Lagrangian formulation of
the constrained minimization,
the mixed finite element method, the mortar domain decomposition method
and the discretization of elliptic and parabolic interface problems.
The two variable parameters can be updated at each iteration,
requiring no a priori estimates on the spectrum
of two preconditioned subsystems involved.
The convergence and convergence rate of the algorithm
are analysed. Both symmetric and nonsymmetric
saddle-point systems are discussed, and numerical experiments
are presented to demonstrate the robustness and effectiveness of the
algorithm.
\end{abstract}

\section{Introduction}

The aim of the current work is to develop
an inexact preconditioned Uzawa algorithm  for the
generalized saddle-point problem of the form
\begin{equation} \label{saddle}
\left( \begin{array}{cr} A & B \\ \\ B^t & -D \end{array} \right)
\left( \begin{array}{c} x \\ \\ y \end{array} \right) =\left(
\begin{array}{c} f \\ \\ g \end{array}\right),
\end{equation}
where $A$ is an $n\times n$
symmetric and positive definite matrix,
$B$ is an $n \times m$ matrix, and
  $D$ is an $m\times m$
symmetric positive semi-definite matrix. We shall assume that the
Schur complement matrix
$$
S=B^tA^{-1}B+D
$$
associated with the system (\ref{saddle}) is an $m\times m$
symmetric and positive definite matrix, which ensures the unique
solvability of system (\ref{saddle}). System \eqref{saddle} arises
from many areas of computational sciences and engineerings, such as
the constrained optimization, the mixed finite element formulation
for the second order elliptic equation, the linear elasticity
problem, as well as elliptic and parabolic interface problems; see
\cite{s6} \cite{chenduzou} \cite{s14} \cite{GloLeta89} \cite{s15}
and Section~\ref{sec:appls}
for several such applications.

Many numerical methods such as Schur complement reduction methods,
null space methods, penalty methods, multilevel methods, Krylov
subspace methods and preconditioning, are investigated to solve the
saddle point problem \eqref{saddle}, especially for solving the
simplest case of the saddle-point system (\ref{saddle}) when the (2,
2) block $D$ vanishes; see \cite{s3} \cite{s4} \cite{s5}
\cite{GloLeta89} \cite{s17}  \cite{s18}
\cite{BenziGolubLiesen_Acta05} and the references therein. In
particular, the inexact preconditioned Uzawa-type algorithms have
attracted wide attention; see \cite{s3} \cite{s4} \cite{s5}
\cite{s13} \cite{huzou01} \cite{huzou} \cite{huzou06} \cite{s17},
and the references therein. These inexact Uzawa-type algorithms have
an important feature that they preserve the minimal memory
requirement and do not need actions of the inverse matrix $A^{-1}$.
On the contrary, few studies on the convergence analysis of
inexact preconditioned Uzawa iterative methods can be found in
the literature
for the generalized saddle-point system (\ref{saddle})
where a general block $D$ is present.
This work intends to make some initial efforts to fill in the gap.

Suppose that $\hat{A}$ and $\hat{S}$ are two symmetric and positive
definite matrices, and act as the preconditioners for $A$ and $S$,
respectively. We shall be interested in the following inexact
preconditioned Uzawa method for solving the system (\ref{saddle}).
%

\begin{equation}\label{eq:alg2}
\begin{array}{l}
x_{i+1}=x_{i}+\omega_i\,\hat{A}^{-1}(f-Ax_{i}-By_{i})\,,
\\
y_{i+1}=y_{i}+\tau_i\hat{S}^{-1}(B^tx_{i+1}-Dy_{i}-g)\,,
\end{array}
\end{equation}
where $\omega_i$ and $\tau_i$ are two relaxation parameters to be
determined at each iteration. Equivalently, the system
(\ref{eq:alg2}) can be written as
\begin{equation*}\label{eq:alg1}
\begin{array}{l}
\ds \hat{A}\frac{x_{i+1}-x_{i}}{\omega_i}+Ax_{i}+By_{i}=f\,,
\\
\ds -\hat{S}\frac{y_{i+1}-y_{i}}{\tau_i}+B^tx_{i+1}-Dy_{i}=g.
\end{array}
\end{equation*}
We shall often need the approximate Schur complement of $S$, namely
%
%
$$
H=B^t\hat{A}^{-1}B+D.
$$

The inexact preconditioned Uzawa method (\ref{eq:alg2}) with two
variable relaxation parameters was first proposed and analysed in
\cite{huzou01} for the simple case of $D=0$, and different variants
of the algorithm were further studied in \cite{huzou}
\cite{huzou06}. The original idea of introducing the variable
parameters $\omega_i$ and $\tau_i$ was to ensure that the resulting
inexact preconditioned Uzawa algorithms always converge for any
available symmetric and positive definite preconditioners $\hat{A}$
and $\hat{S}$, and converge nicely when effective preconditioners
are available. Nearly all other existing preconditioned Uzawa
algorithms do not adopt self-updating relaxation parameters, and
converge only under some proper scalings of the preconditioners
$\hat A$ and $\hat S$.

The choice of the relaxation parameters $\omega_i$ and $\tau_i$ in
(\ref{eq:alg2}) is not straightforward. They should be easily
updated at each iteration and their evaluations should be less
expensive.  The usual choices of parameters by minimizing the errors
$x-x_i$ and $y-y_i$ in certain norms do not work since the
evaluation of the resulting parameters always involve the action of
$A^{-1}$; see \cite{huzou01} for details.

Next, we follow \cite{huzou01} to work out an effective way to
evaluate the two relaxation parameters $\omega_i$ and $\tau_i$ in
(\ref{eq:alg2}). To do so, we consider the two residuals associated
with the $i$-th iteration:
\begin{equation}  \label{res}
f_i=f-(Ax_i+By_i),\quad g_i=B^t x_{i+1}-Dy_i-g .
\end{equation}
Then we may determine the parameter $\omega_i$ by minimizing
$$
|\omega_i\hat{A}^{-1}f_i-A^{-1}f_i |_A^2\,,
$$
which yields
\begin{equation}\label{eq:omega1}
\omega_i=\frac{\langle f_i,r_i \rangle}{ \langle Ar_i,r_i \rangle},
\end{equation}
where $r_i=\hat{A}^{-1}f_i$\, , and $\langle \cdot,\cdot \rangle$
stands for the inner product of two vectors in Euclidean space. The
parameter $\tau_i$ can be determined by minimizing
$$
|\tau_i \hat{S}_i^{-1} g_i-H^{-1} g_i |_H^2,
$$
which gives a prototype choice  \bb\l{eq:tau1}
\hat{\tau}_i=\frac{\langle g_i,s_i \rangle} {\langle Hs_i,s_i
\rangle} \ee with $s_i=\hat{S}^{-1}g_i$. But as we shall see, such
choice of $\tau_i$ may not guarantee the convergence of
Algorithm~\ref{algo:LinInexUzawaVariableRelax}. We need a damping
factor $\theta_i$ for the parameter $\hat\tau_i$ in (\ref{eq:tau1}),
and will take $\tau_i$ in (\ref{eq:alg2}) as
\begin{equation} \label{tau2}
\tau_i=\theta_i\hat\tau_i=\theta_i\frac{\langle g_i,s_i \rangle}
{\langle Hs_i,s_i \rangle}\,.
\end{equation}
The algorithm can be summarized as follows.
\begin{algorithm}
\caption{Linear inexact Uzawa algorithm with variable relaxation. }
\label{algo:LinInexUzawaVariableRelax}
\begin{algorithmic}
\STATE{
\begin{enumerate}
\item Compute $f_i = f-(Ax_{i}+By_{i})$, $r_i=\hat{A}^{-1} f_i $, and $\omega_i=\frac{\langle f_i,r_i \rangle}{\langle Ar_i,r_i \rangle}$; 
\item Update $x_{i+1}=x_{i}+\omega_i\,r_i$;
\item Compute $g_i=B^t x_{i+1} - D y_i -g$, $s_i=\hat{S}^{-1}g_i$,
and $\tau_i=\theta_i \frac{\langle g_i,s_i \rangle}{\langle H
s_i,s_i \rangle}$;
%
\item Update $y_{i+1}=y_{i}+\tau_i s_i$.
\end{enumerate}
}
\end{algorithmic}
\end{algorithm}

Algorithm \ref{algo:LinInexUzawaVariableRelax} was analyzed in
\cite{huzou01} for the simplest case of saddle-point problem
\eqref{saddle} when the (2,2) block $D$ vanishes. Unfortunately the
convergence and convergence rate of Algorithm
\ref{algo:LinInexUzawaVariableRelax} were established still under
some appropriate scaling of preconditioner $\hat A$ for $A$, i.e.,
the smallest eigenvalue of the preconditioned system $\hat A^{-1}A$
is larger than one, although no any appropriate scaling of
preconditioner $\hat S$ for Schur Complement $S$ was needed. In this
work we shall extend the analysis in \cite{huzou01} to the more
general and challenging indefinite system \eqref{saddle}, where the
block $D$ is present. As it will be seen, such an extension is
highly nontrivial for a general block $D$. We need to make essential
modifications of the major analysis techniques in \cite{huzou01} and
introduce several crucial new techniques in order to succeed in
analyzing the convergence and convergence rate of Algorithm
\ref{algo:LinInexUzawaVariableRelax} for general $D\ne 0$. It is
important to remark that for the case of $D=0$, our subsequent
analysis will improve the convergence results, relax the convergence
conditions in \cite{huzou01} and provide instructive information on
the selection of the damping parameter $\theta_i$ to ensure the
convergence. Unlike in \cite{huzou01}, we will not assume
appropriate scalings of two preconditioners $\hat A$ and $\hat S$
for the convergence of Algorithm
\ref{algo:LinInexUzawaVariableRelax}.

We will also generalize Algorithm
\ref{algo:LinInexUzawaVariableRelax} to the cases when the action of
preconditioner $\hat{A}$ or $\hat{S}$ is replaced by a nonlinear
iterative solver. This is more practical and important for some
applications where effective preconditioners are not available.
The proposed algorithm is also analyzed and tested numerically when $A$
in \eqref{saddle} is nonsymmetric. No such analysis is available
in the literature when inexact preconditioners are used.

For the sake of clarity, we list the main notations used later.
\begin{longtable}{lp{4.0in}}\hline
\multicolumn{2}{ l }{\bf Some notations and definitions} \\
$S$          & $S=B^t A^{-1} B + D$, and $S = R R^t$ with $R$ being nonsingular \\
$\hat{S}$          & an spd approximation of $S$  \\
$H$             & $H = B^t \hat{A}^{-1} B + D$, where $\hat{A}$ is an approximation of $A$ \\
$\theta_i$             &  damping factor for  the parameter $\hat{\tau}_i$  \\
$\alpha$             &  $\alpha = (\kappa_1 -1 ) / (\kappa_1 + 1)$, where $\kappa_1 = \text{cond}(\hat{A}^{-1} A)$ \\
$\beta$             &  $\beta = (\kappa_2 -1 ) / (\kappa_2 + 1)$, where $\kappa_2 = \text{cond}(\hat{S}^{-1} H)$ \\
$c_0$         & the largest eigenvalue of $D^{-1} B^t \hat{A}^{-1} B$ (see Lemma \ref{lem:eig_estimate1}) \\
$c_1$           & the constant defined in \eqref{weigh} or \eqref{eq:weigh1} \\
$Q_i$           & defined by $Q_i^{-1} = \theta_i G_i^{-1}$ in \eqref{eq:definitionQi}, and $G_i$ is given in  Lemma \ref{lem:eig_estimate1}  \\
$\lambda, \lambda_0$           & $\lambda \hat{A} \leq A \leq \lambda_0 \hat{A}$;  see \eqref{(1)} \\
$\delta_1, \delta_2$           & $ \delta_1 = (\lambda_0 + c_0)/ [\lambda_0 ( 1+c_0)] $, $ \delta_2 = (\lambda + c_0)/ [\lambda ( 1+c_0)]$ (see Lemma \ref{lem:eig_estimate1}) \\
$\omega, \Omega$     & $\omega |z|^2 \leq |W z|^2 \leq \Omega |z|^2, \forall z$; see \eqref{eq:w_bound1} \\
\hline
\end{longtable}

\section{Basic formulation}
\setcounter{equation}{0}
We shall often use the condition
numbers of the two preconditioned systems
$$\begin{array}{l}
\kappa_1=\m{cond}(\hat{A}^{-1}A),\qquad
\kappa_2=\m{cond}(\hat{S}^{-1}H)
\end{array} $$
and the following two convergence-rate related constants
$$\begin{array}{l}
\alpha=\frac{\kappa_1-1}{\kappa_1+1}\,, \qquad
\beta=\frac{\kappa_2-1}{\kappa_2+1}.
\end{array} $$
For any two symmetric and semi-positive definite matrices $C_1$ and
$C_2$ of order $m$ satisfying
$$
\langle C_1\phi, \phi \rangle \le \langle C_2\phi, \phi \rangle ,
\quad \forall\,\phi\in \boldmath \mathbb{R}^m\,,
$$
we will simply write
$$
C_1\le C_2\,.
$$

\subsection{When the (2, 2) block $D$ vanishes}

The convergence of Algorithm \ref{algo:LinInexUzawaVariableRelax}
was analyzed in \cite{huzou01} for the saddle-point system
(\ref{saddle}) with $D=0$ under the condition that preconditioner
$\hat{A}$  for $A$ is appropriately scaled such that
\begin{equation}\label{eq:huzoucond}
\hat{A} \le A \le \hat\lambda_0 \hat{A}
\end{equation}
for some constant $\hat\lambda_0\ge 1$. Next, we will demonstrate
the convergence of Algorithm \ref{algo:LinInexUzawaVariableRelax}
without the condition (\ref{eq:huzoucond}).  In fact, since $\hat A$
is a general preconditioner for $A$, there are always two positive
constants $\lambda$ and $\lambda_0$ such that $\lambda\le 1
\le\lambda_0$ and
\begin{equation}  \label{(1)}
\lambda\,\hat{A} \le A  \le \lambda_0\,\hat{A}.
\end{equation}
Noting that (\ref{(1)}) is not an actual assumption
since it is always true.
Now we let
\begin{equation}  \label{eq:hatA}
\tilde{A}=\lambda\,\hat{A}\,, \quad
\tilde{\lambda}_0={\lambda_0}/{\lambda}.
\end{equation}
Then we can rewrite (\ref{(1)}) as
\begin{equation} \label{(2)}
\tilde{A} \le A \le \tilde{\lambda}_0\, \tilde{A}\,.
\end{equation}
In terms of this newly introduced $\tilde A$, one may express
Algorithm \ref{algo:LinInexUzawaVariableRelax} as follows (noting
that $D=0$):
\begin{equation*}\label{eq:alg3}
\begin{array}{l}
x_{i+1}=x_{i}+\tilde{\omega}_i\,\tilde{A}^{-1}(f-Ax_{i}-By_{i})\,,
\\ \\
y_{i+1}=y_{i}+\tilde{\tau}_i\tilde{S}^{-1}(B^tx_{i+1}-g),
\end{array}
\end{equation*}
where $\tilde{S}=\lambda\, \hat{S}$ and the damping parameters $\tilde \omega_i$ and $\tilde \tau_i$
are given by
%
%
\begin{equation}\label{eq:omega3}
\tilde{\omega}_i=\lambda\,\omega_i, \quad
\tilde{\tau}_i=\lambda \,\tau_i\,.
\end{equation}

Let us introduce a parameter
\begin{equation} \label{alph}
\alpha_i= \frac{|(I-\tilde{\omega}_iA^{\frac{1}{2}}\tilde{A}^{-1}
A^{\frac{1}{2}})A^{-\frac{1}{2}}f_i|}{|A^{-\frac{1}{2}}f_i|} =
\frac{|(I-{\omega}_iA^{\frac{1}{2}}\hat{A}^{-1}
A^{\frac{1}{2}})A^{-\frac{1}{2}}f_i|}{|A^{-\frac{1}{2}}f_i|}\,,
\end{equation}
where the residual $f_i$ is defined by \eqref{res}. Then we can show

\begin{lem}\label{lem:first} For the parameters $\tilde{\lambda}_0$, $\tilde{\omega}_i$
and $\alpha_i$ defined respectively in \eqref{eq:hatA},
\eqref{eq:omega3} and \eqref{alph}, it holds that
$$
\tilde{\lambda}_0^{-1}\le \tilde{\omega}_i\le 1-\alpha_i^2,\quad 0
\le \alpha_i \le \alpha.
$$
\end{lem}

\noindent {Proof}. First note that if we follow the same way as we
get $\omega_i$ in (\ref{eq:omega1}) when $\hat A$ is replaced by
$\tilde A$, we derive a new parameter $\tilde \omega_i$, which is
exactly the one given by $\tilde{\omega}_i=\lambda\,\omega_i$ in
(\ref{eq:omega3}). Then following the proof of Lemma~3.1 in
\cite{huzou01} by means of the relations (\ref{(2)}) and the fact
that cond($\tilde A^{-1}A$)=cond($\hat A^{-1}A$), we can derive the
desired estimates. $\square$

\ss Intuitively it is easy to understand that Algorithm
\ref{algo:LinInexUzawaVariableRelax} may not converge for an
arbitrary damping parameter $\theta_i$ in (\ref{tau2}). Following
the convergence analysis in \cite{huzou01} and using
Lemma~\ref{lem:first}, we have the following convergence.
\begin{lem}\label{lem:D1}
For any damping parameter $\theta_i$ in \eqref{tau2} satisfying
\begin{equation}\label{eq:cond_theta1}
\theta_i(1+\beta)\le 1-\alpha_i\,,
\end{equation}
Algorithm \ref{algo:LinInexUzawaVariableRelax} converges.
\end{lem}

\begin{rem} The formula \eqref{eq:cond_theta1} gives a range to choose
the damping parameter $\theta_i$ for the convergence of Algorithm
\ref{algo:LinInexUzawaVariableRelax}, but it does not provide the
best choice of $\theta_i$.

If one can estimate the lower bound $\lambda$ in \eqref{(1)}, say
$\hat \lambda$ is such an estimate.   That is, $\hat{\lambda}
\omega_i \leq \lambda \omega_i = \tilde{\omega}_i$. Therefore, $1-
\sqrt{1-\hat{\lambda}\,\omega_i} \le 1-\sqrt{1-\tilde{\omega}_i} \le
1-\alpha_i$.   Then we can have a more explicit range for $\theta_i$
to ensure the convergence:
\begin{equation} \label{Th1}
\theta_i\le \frac{1-\sqrt{1-\hat{\lambda}\,\omega_i}}{2}\,.
\end{equation}
In fact, it follows from \eqref{Th1} that $ \theta_i(1+\beta) \leq 2
\theta_i \le 1- \sqrt{1-\hat{\lambda}\,\omega_i} 
\le 1-\alpha_i\,, $  hence Algorithm
\ref{algo:LinInexUzawaVariableRelax} converges by
Lemma~\ref{lem:D1}. A lower bound estimate $\hat{\lambda}$ may be
obtained by knowing a upper bound $\hat{\kappa}_1$ of the condition
number $\ds\kappa_1={\lambda_0}/{\lambda}$ and a lower bound
$\hat{\lambda}_0$ of $\lambda_0$ and letting
$$
\hat{\lambda}=\frac{\hat{\lambda}_0}{\hat{\kappa}_1}.
$$
The lower bound $\hat \lambda_0$ can be easily evaluated, e.g.,
using the power method for $\hat{A}^{-1}A$.
\end{rem}

\subsection{When the (2, 2) block $D$ is present}

In this subsection we will study the convergence of Algorithm
\ref{algo:LinInexUzawaVariableRelax}  when the block matrix $D$ is
present. As we will see, the convergence of Algorithm
\ref{algo:LinInexUzawaVariableRelax} is much more complicated than
the case with $D=0$, and it is essential to scale preconditioner
$\hat{A}$ for $D\ne 0$ since the convergence of Algorithm
\ref{algo:LinInexUzawaVariableRelax} depends strongly on the
relative scale of $B^t\hat{A}^{-1}B$ with respect to $D$ in the
approximated Schur complement $H=B^t\hat{A}^{-1}B+D$. The following
lemma illustrates this fact in terms of the eigenvalues of the
preconditioned Schur complement and is essential to the subsequent
convergence analysis.

\begin{lem}\label{lem:eig_estimate1}
Let
$$
\beta_i=\frac{|(I-\hat{\tau}_iH^{\frac{1}{2}}\hat{S}^{-1}
H^{\frac{1}{2}})H^{-\frac{1}{2}}g_i|}{|H^{-\frac{1}{2}}g_i|}\,,
$$
where $g_i$ is defined in \eqref{res}, then it holds
$$
1-\beta_i^2=\hat{\tau}_i\frac{\langle g_i,\hat{S}^{-1}g_i \rangle
}{\langle g_i,H^{-1}g_i \rangle } \quad\mbox{and}\quad 0 \le \beta_i
\le \beta.
$$
And there exists a symmetric and positive definite matrix $G_i$ such
that $G_i^{-1}g_i=\hat\tau_i\hat S^{-1}g_i$ and all the eigenvalues
of the preconditioned matrix $ S^{\frac 12 }G_i^{-1} S^{\frac 12}$
(or $ R^t G_i^{-1} R$, where $S = R R^t$) lie in the interval
$$
[  (1-\beta_i) \delta_1, ~  (1+\beta_i) \delta_2 ],$$ where $
\delta_1=\frac{\lambda_0+c_0}{\lambda_0(1+c_0)}$,
$\delta_2=\frac{\lambda +c_0}{\lambda(1+c_0)}$, and $c_0$ is the
largest eigenvalue of $D^{-1}B^t\hat{A}^{-1}B$.
\end{lem}

\noindent{Proof}. By the definition of $\hat{\tau}_i$, we can write
\begin{eqnarray}
|\hat{\tau}_i\hat{S}^{-1}g_i-H^{-1}g_i|^2_H &=&
\hat\tau_i^2|\hat{S}^{-1}g_i|^2_H-2\hat{\tau}_i \langle
g_i,\hat{S}^{-1}g_i \rangle
+|H^{-1}g_i|_H^2\nonumber\\
&=&\Big(1-\hat{\tau}_i\frac{ \langle g_i,\hat{S}^{-1}g_i \rangle }{
\langle g_i,H^{-1}g_i \rangle } \Big)\,|H^{-1}g_i|^2_H\,,\l{eq:bi1}
\end{eqnarray}
and
\begin{eqnarray*}
\ds |(I-\hat{\tau}_iH^{\frac{1}{2}}\hat{S}^{-1}
H^{\frac{1}{2}})H^{-\frac{1}{2}}g_i|^2 &=&
|H^{-\frac{1}{2}}g_i|^2-2\hat{\tau}_i \langle g_i,\hat{S}^{-1}g_i
\rangle +\hat{\tau}_i^2 \langle H\hat{S}^{-1}g_i,\hat{S}^{-1}g_i
\rangle
\\
&=&\Big(1-\hat{\tau}_i\frac{ \langle g_i,\hat{S}^{-1}g_i \rangle }{
\langle g_i,H^{-1}g_i \rangle } \Big)\,|H^{-1}g_i|^2_H\,.
\end{eqnarray*}
Thus it follows from the definition of $\beta_i$ and
(\ref{eq:bi1}) that
$$
\beta_i^2=1-\hat{\tau}_i\frac{ \langle g_i,\hat{S}^{-1}g_i \rangle
}{ \langle g_i,H^{-1}g_i \rangle }
$$
and
\begin{equation}\label{eq:betai}
|\hat{\tau}_i \hat{S}^{-1}g_i-H^{-1}g_i|_H = 
\beta_i\,|H^{-1}g_i|_H.
\end{equation}
It is shown in the proof of Lemma\,3.2  in \cite{huzou01} by using
the Kantorovich inequality that
$$
\hat{\tau}_i\frac{ \langle g_i,\hat{S}^{-1}g_i \rangle }{ \langle
g_i,H^{-1}g_i \rangle } \ge 1-\beta^2
$$
and thus we have $\beta_i\le \beta$.
On the other hand, the estimate (\ref{eq:betai})
implies the existence of a symmetric and positive definite matrix $G_i$ such
that (cf.~\cite{s3})
$$
G_i^{-1}g_i=\hat{\tau}_i\hat{S}^{-1}g_i
$$
and
$$
|I-H^{\frac{1}{2}}G_i^{-1}H^{\frac{1}{2}}|\le \beta_i.
$$
Let $\mu>0$ be an eigenvalue of $S^{\frac{1}{2}}
G_i^{-1}S^{\frac{1}{2}}$, then there exists a vector $\phi$ such
that
$$
\langle S\phi,\phi \rangle =\mu \langle G_i\phi,\phi \rangle .
$$
It is easy to see that for any $\phi\in \boldmath \mathbb{R}^m$,
$$
 \langle B^t {A}^{-1}B\phi,\phi  \rangle =
 \langle (\hat{A}^{\frac{1}{2}}A^{-1}\hat{A}^{\frac{1}{2}})\hat{A}^{-\frac{1}{2}}B\phi,
\hat{A}^{-\frac{1}{2}}B\phi  \rangle .
$$
But we know from the assumption \eqref{(1)} that
$$
\frac{1}{\lambda_0}I\le
\hat{A}^{\frac{1}{2}}A^{-1}\hat{A}^{\frac{1}{2}} \le
\frac{1}{\lambda}I\,,
$$
which leads to
\begin{equation} \label{est}
 \langle (\frac{1}{\lambda_0}B^t\hat{A}^{-1}B+D)\phi,\phi \rangle  \le  \langle (B^t
A^{-1}B+D)\phi,\phi \rangle  \le
 \langle (\frac{1}{\lambda}B^t\hat{A}^{-1}B+D)\phi,\phi \rangle .
\end{equation}
Let $\gamma\in (0, 1)$ be a constant to be determined such that for
any $\phi\in \boldmath \mathbb{R}^m$,
$$
 \langle (B^t\hat{A}^{-1}B+\lambda D)\phi, \phi \rangle  \le \gamma
 \langle H\phi, \phi  \rangle =\gamma  \langle B^t\hat{A}^{-1}B\phi, \phi \rangle  +\gamma  \langle D\phi,
\phi \rangle ,
$$
which implies
$$
 \langle B^t\hat{A}^{-1}B\phi, \phi \rangle  \le \frac{\gamma-\lambda}{1-\gamma}
 \langle D\phi, \phi \rangle .
$$
As $c_0$ is the largest eigenvalue of
$D^{-1}B^t\hat{A}^{-1}B$, we can choose $\gamma$ such that
${(\gamma-\lambda)}/{(1-\gamma)}=c_0$,
that gives $\gamma = (\lambda+c_0)/(1+c_0)$. Hence we know
$$
 \langle (B^t\hat{A}^{-1}B+\lambda D)\phi, \phi \rangle  \le \frac{\lambda+c_0}{1+c_0}
 \langle H\phi, \phi \rangle \,.
$$
Similarly we can derive
$$
 \langle (B^t\hat{A}^{-1}B+\lambda_0 D)\phi, \phi \rangle  \ge  \frac{\lambda_0+c_0}{1+c_0}
 \langle H\phi, \phi \rangle \,.
$$
Using the above two estimates we deduce from (\ref{est}) that
$$
\frac{1}{\lambda_0}\frac{\lambda_0+c_0}{1+c_0} \langle H\phi,\phi
\rangle \le  \langle S \phi, \phi \rangle  =  \mu\, \langle
G_i\phi,\phi \rangle  \le \frac{1}{\lambda}\frac{\lambda+c_0}{1+c_0}
 \langle H\phi,\phi \rangle .
$$
Since $|I-H^{\frac{1}{2}}G_i^{-1}H^{\frac{1}{2}}|\le \beta_i$,
$$
\frac{1-\beta_i}{\lambda_0}\frac{\lambda_0+c_0}{1+c_0}  \langle
G_i\phi,\phi \rangle  \le \mu\, \langle G_i\phi,\phi \rangle  \le
\frac{1+\beta_i}{\lambda}\frac{\lambda+c_0}{1+c_0} \langle
G_i\phi,\phi \rangle .
$$
which implies the claimed eigenvalue bounds. $\square$ \vspace{2mm}

\begin{rem}
We give some comments on the constant $c_0$ introduced in
Lemma~\ref{lem:eig_estimate1}. If $D=0$, then $c_0=\infty$ and the
estimate in Lemma~\ref{lem:eig_estimate1} coincides with the one in
\cite{huzou01}. Noting that $\lambda A^{{-1}} \le \hat{A}^{-1}\le
\lambda_0 A^{-1}$, $c_0$ is bounded above and below by the
eigenvalues of $D^{-1}B^tA^{-1}B$. Thus if $D$ dominates, then $c_0$
is very small.
\end{rem}

\begin{rem}
Lemma~\ref{lem:eig_estimate1} shows that unless $D$ dominates
$B^t\hat{A}^{-1}B$ (i.e., $c_0$ is small) the eigenvalues  of the
preconditioned matrix $\theta_i
S^{\frac{1}{2}}G_i^{-1}S^{\frac{1}{2}}$ are  sensitive to the
scaling of preconditioner $\hat{A}$. Thus, we may assume $\lambda=1$
as in \eqref{(2)}, i.e.,
$$
\hat{A} \le A \le \lambda_0 \hat{A},
$$
which may be achieved by scaling $\hat{A}$ by $\lambda$; see
\eqref{eq:hatA}. Then we may simply choose the damping parameter
$\theta_i$ in \eqref{tau2} such that
\begin{equation*}
\theta_i\le \frac{1-\sqrt{1-\omega_i}}{2}
\end{equation*}
for the convergence of Algorithm
\ref{algo:LinInexUzawaVariableRelax}.
In general we may choose the damping parameter $
\theta_i={M}/{\kappa_1} $, where the constant $M$ should be selected
for guaranteeing the convergence of Algorithm
\ref{algo:LinInexUzawaVariableRelax} and achieving an appropriate
convergence rate; we refer to the further discussions in the next
section.
\end{rem}

\section{Convergence analysis}

Now, we are ready to analyze the convergence of Algorithm
\ref{algo:LinInexUzawaVariableRelax}. Let us introduce the errors
$$
e^x_i=x-x_i,\quad e^y_i=y-y_i\,.
$$
Then the residuals $f_i$ and $g_i$ can be expressed as \bb \l{eq:fi}
f_i=Ae^x_i+Be^y_i\,, \quad g_i=-B^te_{i+1}^x+De_i^y\,. \ee Using the
definition of $f_i$ and the iteration (\ref{eq:alg2}) for updating
$x_i$, we can write \bb \l{eq:ei1} A^{\frac{1}{2}}e^x_{i+1}=
A^{\frac{1}{2}}(e^x_i-\omega_i\hat{A}^{-1}f_i) =(I-\omega_i
A^{\frac{1}{2}}\hat{A}^{-1}A^{\frac{1}{2}})A^{-\frac{1}{2}}f_i
-A^{-\frac{1}{2}}Be^y_i\,. \ee On the other hand, using the
iteration (\ref{eq:alg2}) for updating $y_i$, the definition of
$g_i$, the formula (\ref{eq:ei1}) and the matrix $G_i$ introduced in
Lemma~\ref{lem:eig_estimate1} we derive \bb \l{eq:long1b}
\begin{array}{l}
\ds e^y_{i+1}=e^y_i-\tau_i \hat S^{-1}g_i
=e^y_i-Q_i^{-1}g_i=e_i^y+Q_i^{-1}(B^te^x_{i+1}-De^y_i)
\\ \\
\ds \quad =e^y_i+Q_i^{-1}B^tA^{-\frac{1}{2}} \left((I-\omega_i
A^{\frac{1}{2}}\hat{A}^{-1}
A^{\frac{1}{2}})A^{-\frac{1}{2}}f_i-A^{-\frac{1}{2}}Be^y_i \right)
-Q_i^{-1}De^y_i
\\ \\
\ds\quad =Q_i^{-1}B^tA^{-\frac{1}{2}}(I-\omega_i
A^{\frac{1}{2}}\hat{A}^{-1}A^{\frac{1}{2}})A^{-\frac{1}{2}}f_i
+(I-Q_i^{-1}S)e^y_i\,,
\end{array}
\ee
where
\begin{equation}\label{eq:definitionQi}
Q_i^{-1}\equiv {\tau}_i\hat{S}^{-1}
=\theta_i\hat{\tau}_i\hat{S}^{-1} =\theta_i G_i^{-1} \,.
\end{equation}
Now it follows from (\ref{eq:fi}), (\ref{eq:ei1}) and
(\ref{eq:long1b}) that
\bb \l{eq:fi+1b}
\begin{array}{l}
A^{-\frac{1}{2}}f_{i+1}=A^{\frac{1}{2}}e^x_{i+1}+A^{-\frac{1}{2}}Be^y_{i+1}
\\ \\
\quad =(I+A^{-\frac{1}{2}}BQ_i^{-1}B^tA^{-\frac{1}{2}})
(I-\omega_i
A^{\frac{1}{2}}\hat{A}^{-1}A^{\frac{1}{2}})A^{-\frac{1}{2}}f_i
-A^{-\frac{1}{2}}BQ_i^{-1}Se^y_i\,.
\end{array} \ee

Consider the singular value decomposition of the matrix $B^tA^{-\frac{1}{2}}$,
$$
B^tA^{-\frac{1}{2}}=U\Sigma V^t,\quad \Sigma=[\Sigma_0, \; 0],
$$
where $U$ is an $m\times m$ orthogonal matrix, $V$ is an
$n\times n$
orthogonal matrix, and $\Sigma_0$ is an $m\times m$
diagonal matrix with its diagonal entries being the singular values of
$B^tA^{-\frac{1}{2}}$, and $S=RR^t$, where $R$ is a non-singular $m\times m$ matrix.
Set
$$
E^{(1)}_i=\sqrt{\alpha}V^tA^{-\frac{1}{2}}f_i,\quad E^{(2)}_i=R^t
e^y_i\,,
$$
then we can write by using (\ref{eq:long1b}) and (\ref{eq:fi+1b}) that
$$
\left(\begin{array}{c} E^{(1)}_{i+1} \\ \\ E^{(2)}_{i+1} \end{array}
\right)= \left(\begin{array}{cc} \alpha(I+\Sigma^t
U^tQ_i^{-1}U\Sigma) & \sqrt{\alpha}
\Sigma^t U^tQ_i^{-1} R  \\ \\
\sqrt{\alpha} R^t Q_i^{-1}U\Sigma &
-(I- R^t Q_i^{-1} R )
\end{array} \right) \left( \begin{array}{c} \frac{1}{\alpha}V^t
(I-\omega_i A^{\frac{1}{2}}\hat{A}^{-1}A^{\frac{1}{2}})VE^{(1)}_i \\
\\ -E^{(2)}_i
\end{array}\right)\, .
$$
One can easily verify by noting $\alpha_i\le \alpha$
that
$$
 |\frac{1}{\alpha}V^t(I-\omega_i
A^{\frac{1}{2}}\hat{A}^{-1}A^{\frac{1}{2}})VE^{(1)}_i |^2
=\frac{\alpha_i^2}{\alpha^2} |E^{(1)}_i |^2 \le |E^{(1)}_i |^2\,,
$$
which, along with the relations
$$
(I+\Sigma^t U^t Q_i^{-1}U \Sigma )=\left(\begin{array}{cc}
I+\Sigma_0^t U^tQ_i^{-1}U \Sigma_0  & 0 \\ \\ 0 & I
\end{array}\right),\quad
\Sigma^t U^t Q_i^{-1} R
=\left(\begin{array}{c} \Sigma_0^t U^tQ_i^{-1} R \\
\\
0 \end{array}\right),
$$
enables us to reduce the estimate of components $E^{(1)}_{i+1}$ and
$E^{(2)}_{i+1}$ to the spectral estimate of the following symmetric
matrix
\begin{equation}\label{eq:Fi}
F_i=\left( \begin{array}{cc} \alpha\,(I+\Sigma_0^t U^t Q_i^{-1} U
\Sigma_0) &
\sqrt{\alpha}\, \Sigma_0^t U^tQ_i^{-1} R \\ \\
\sqrt{\alpha}\, R^t Q_i^{-1} U \Sigma_0 &
-(I-R^tQ_i^{-1}R)\end{array} \right)\,.
\end{equation}
So if all the eigenvalues of $F_i$ are bounded by $\rho= ||F_i ||<1$
in their magnitude, then
\begin{equation}\label{eq:rate1}
|E^{(1)}_{i+1}|^2+|E^{(2)}_{i+1}|^2\le
\rho^2\,(\frac{\alpha_i^2}{\alpha^2}|E^{(1)}_{i}|^2+|E^{(2)}_{i}|^2),
\end{equation}
 and the convergence of Algorithm
\ref{algo:LinInexUzawaVariableRelax} can be ensured. In following we
will have an estimate on $\rho$. We first examine the convergence
the algorithm, and then further estimate the convergence rate.

For the spectral estimate of $F_i$, we introduce a parameter $c_1$
satisfying
\begin{equation} \label{weigh}
c_1\, R^t Q_i^{-1}R \ge \Sigma_0^tU^tQ_i^{-1}U\Sigma_0,
\end{equation}
where $c_1$ measures the magnitude of $B^tA^{-1}B$ relatively to the one of
$D$ in an appropriately weighted sense.
If we let
$$
W=Q_i^{-\frac{1}{2}} R, \quad T=Q_i^{-\frac{1}{2}}
U\Sigma_0\,,
$$
then (\ref{weigh}) is equivalent to the following inequality
\bb\label{eq:weigh1}  |T u |^2 \le c_1\, |W u |^2\, \quad \forall\,
u\in \mathbb{R}^m\,. \ee Using the parameter $c_1$, we have the
following result.

\begin{thm}\label{thm:spectral_f1} 
If the damping
parameter $\theta_i$ satisfies
$$
\theta_i(1+\beta)\delta_2 <
\frac{2(1-\alpha)}{1-\alpha + 2c_1\alpha}\,,
%
%
$$
then $\rho= ||F_i ||<1$, so Algorithm
\ref{algo:LinInexUzawaVariableRelax} converges.
\end{thm}

\noindent {\it Proof}. We estimate the  upper and lower bounds of
all the eigenvalues of matrix $F_i$.

To see the lower bound of $F_i$, we observe that for any $u, v\in
\mathbb{R}^m$ with one of them being non-zero,
$$\begin{array}{l}
\ds \left \langle (F_i+I) \left(\begin{array}{c} u \\ v \end{array}
\right ),
\left(\begin{array}{c} u \\ v \end{array} \right)\right \rangle  \\ \\
\ds\qquad =(\alpha+1)|u|^2 +\alpha\, \langle Tu, Tu \rangle  +
\sqrt{\alpha}( \langle Wv, Tu \rangle  +  \langle Tu, Wv \rangle ) +
\langle Wv, Wv \rangle
\\ \\
\ds\qquad =(\alpha+1)|u|^2+|\sqrt{\alpha}Tu+Wv|^2> 0\,, \end{array} $$
thus all the eigenvalues of $F_i$ are bounded below by $-1$.

For the upper bound, we consider
$$
J=\begin{array}{l} \left \langle (I-F_i) \left ( \begin{array}{c} u \\
v \end{array} \right), \left(\begin{array}{c} u \\ v \end{array}
\right)\right \rangle
\end{array}
=(1-\alpha)\,|u|^2-|\sqrt{\alpha}Tu+Wv|^2
+2|v|^2\,.
$$
Let $\omega$ and $\gamma$ be the spectral bounds of $W$ given by
\begin{equation}\label{eq:w_bound1}
\omega \,|z|^2 \le |W z|^2\le \gamma \,|z|^2\, \quad \forall\,z\in
\mathbb{R}^m\,,
\end{equation}
then using Young's inequality
we know for all
$\delta>0$,
$$
|\sqrt{\alpha}Tu+Wv|^2\le
(1+\delta)\alpha\,|Tu|^2+(1+\frac{1}{\delta})\,|Wv|^2,
$$
hence it follows from (\ref{eq:weigh1}) and (\ref{eq:w_bound1})
that
\begin{equation}\label{eq:cond_J}
J\ge (1-\alpha)\,|u|^2-c_1\alpha(1+\delta) \gamma \,|u|^2 +
2|v|^2-(1+\frac{1}{\delta}) \gamma |v|^2\,.
\end{equation}
For $J>0$, we need the existence of a $\delta>0$
such that
$$
(1-\alpha)-c_1\alpha(1+\delta) \gamma >0\mbox{  and  }
2-(1+\frac{1}{\delta}) \gamma >0\,,
$$
which is equivalent to
$$
\frac{\gamma}{2-\gamma}<\delta<\frac{1-\alpha}{c_1 \alpha \gamma
}-1\,,
$$
and hence
$$
2(1-\alpha)-(1-\alpha + 2c_1\alpha) \gamma >0,
$$
which requires $\gamma$ to satisfy
\begin{equation}\label{eq:cond_W}
0< \gamma <\frac{2(1-\alpha)}{1-\alpha + 2c_1\alpha}\,.
\end{equation}
Clearly if this condition holds, then it follows from
(\ref{eq:cond_J}) that $I-F_i>0$. Thus all the eigenvalues of $F_i$
have the upper bound $1$. But by Lemma\,\ref{lem:eig_estimate1}, we
know
$$
|Wy|^2 \le \theta_i\,(1+\beta)\delta_2\,|y|^2
$$
which leads to the desired result of Theorem~\ref{thm:spectral_f1}
by taking $\gamma =\theta_i\,(1+\beta)\delta_2$. $\square$
\vspace{2mm}

\begin{rem} We may comment on some direct consequences
of Theorem~\ref{thm:spectral_f1} at the extreme cases of $c_1$ close
to $0$ or $1$. It is easy to see that for $0 <  c_1 \le 1$,
${2(1-\alpha)}/{(1-\alpha + 2c_1\alpha)}$ is monotonically
decreasing with respect to $c_1$, implying that
$$
\frac{2(1-\alpha)}{1+\alpha} \le \frac{2(1-\alpha)}{1-\alpha +
2c_1\alpha} < 2.
$$
Then Theorem~\ref{thm:spectral_f1} ensures the convergence of
Algorithm \ref{algo:LinInexUzawaVariableRelax} for any $\theta_i$
satisfying
$$
\theta_i (1+\beta)\delta_2 < \frac{2(1-\alpha)}{1+\alpha}=\frac{2}{\kappa_1}
$$
in the case that $c_1$ is close to $1$, i.e., $D$ is relatively
small compared to $B^tA^{-1}B$ in the sense of \eqref{weigh}.
 In the case that $c_0$ and $c_1$ are close to $0$, i.e., $D$
dominates $B^tA^{-1}B$, we take $\theta_i$ satisfying
$$
\theta_i (1+\beta) < 2
$$
to guarantee the convergence according to
Theorem~\ref{thm:spectral_f1}, or roughly we take $\theta_i \leq 1$.

\end{rem}

{\bf Estimate of convergence rate}. The following of this section is
devoted to estimating the convergence rate of Algorithm
\ref{algo:LinInexUzawaVariableRelax}.  That is, the more precise
size of $\rho = ||F_i||$ in \eqref{eq:rate1}.  Following exactly the
same arguments as the one for the upper bound of $F_i$ in the proof
of Theorem~\ref{thm:spectral_f1}, we can show that $F_i\le\mu\,I$
for some $\mu\in (\alpha,1)$ provided that there exists a $\delta>0$
such that
$$
(\mu-\alpha)-c_1\alpha(1+\delta) \gamma \ge 0\quad\mbox{and}\quad
\mu+1-(1+\frac{1}{\delta})\, \gamma \ge 0,
$$
where $\gamma$ is the spectral bound of $W$ in (\ref{eq:w_bound1}).
This implies
$$
\frac{\gamma}{\mu+1-\gamma}\le \delta \le \frac{\mu-\alpha}{\alpha
c_1 \gamma} -1\,,
$$
or equivalently
\begin{equation}\label{eq:rate2}
0< \gamma \le \frac{(\mu+1)(\mu-\alpha)}{\alpha c_1
(\mu+1)+\mu-\alpha} \equiv \gamma(\mu,\alpha,c_1)\,.
\end{equation}
That is, if $\gamma \le \gamma(\mu,\alpha,c_1)$ then all the
eigenvalues of $F_i$ are bounded above by $\mu$.

To estimate the lower bound of $F_i$, for any
$\tilde{\mu}\in(0,1)$ and $0<\delta<1$ we can derive
$$\begin{array}{l}
J= \left \langle (F_i+\tilde\mu I) \left(\begin{array}{c} u \\ v
\end{array} \right), \left(\begin{array}{c} u \\ v \end{array}
\right)\right \rangle
=(\alpha+\tilde{\mu})\,|u|^2+|\sqrt{\alpha}Tu+W v|^2
+(\tilde{\mu}-1)|v|^2
\\ \\
\ds\qquad \ge
(\tilde{\mu}+\alpha)\,|u|^2+(1-\frac{1}{\delta})\alpha\,|Tu|^2 +
(\tilde{\mu}-1)|v|^2+(1-\delta)|Wv|^2\,.
\end{array} $$
Using (\ref{eq:w_bound1}) and (\ref{eq:weigh1}), we get $|T u|^2 \le
c_1 \gamma |u|^2$, thus
$$
J \ge (\tilde{\mu}+\alpha+(1-\frac{1}{\delta})c_1 \alpha
\gamma)\,|u|^2 + (\tilde{\mu}-1+(1-\delta) \omega )\,|v|^2\,.
$$
This implies $F_i\ge -\tilde{\mu}\,I$ if there exists a $\delta>0$
such that
$$
(\tilde{\mu}+\alpha+(1-\frac{1}{\delta})c_1 \alpha \gamma)\ge
0,\quad \tilde{\mu}-1+(1-\delta) \omega \ge 0,
$$
or equivalently
$$
\frac{c_1 \alpha \gamma}{\tilde{\mu}+\alpha+\alpha c_1
\gamma}\le\delta\le \frac{\omega +\tilde{\mu}-1}{\omega},
$$
which is equivalent to requiring that 
\begin{equation} \label{R2}
\omega (\tilde \mu, \alpha, c_1, \gamma) \equiv
(1-\tilde{\mu})(1+\frac{c_1 \alpha \gamma }{\tilde{\mu}+\alpha}) \le
\omega,
\end{equation}
then all the eigenvalues of $F_i$ are bounded below by $-\tilde
\mu$. By Lemma~\ref{lem:eig_estimate1}, we know that
$$ \theta_i (1-\beta) \delta_1 |y|^2 \leq | W y |^2,$$
so we
can take $\omega = \theta_i (1-\beta) \delta_1$.

Note that in \eqref{R2} $ \omega (\tilde \mu, \alpha, c_1, \gamma)
\leq (1-\tilde{\mu}) (1+c_1\,\frac{\alpha}{\alpha+\beta}
\gamma(\mu,\alpha,c_1) ) $ for $\tilde{\mu}\in [\beta,1)$. Using
Lemma~\ref{lem:eig_estimate1}, we derive immediately from
(\ref{eq:rate2})--\eqref{R2} the following results.
\begin{thm}\label{thm:rates1} For any $\mu\in (\alpha,1)$, if $\theta_i$ satisfies
$$
\theta_i \le \frac{1}{\delta_2(1+\beta)} \gamma(\mu,\alpha,c_1)\,,
$$
and $\tilde{\mu}\in [\beta,1)$ satisfies
$$
1-\tilde{\mu}
\le \frac{\delta_1\,(1-\beta)}{1+c_1\,\frac{\alpha}{\alpha+\beta}
\gamma(\mu,\alpha,c_1)} \theta_i,
$$
then we have
$$
-\tilde{\mu}I \le F_i\le \mu\,I\, ,
$$
 and the convergence rate $\rho = \max\{\mu, \tilde{\mu} \}$. 
\end{thm}


{\bf Rate estimates at extreme cases}. We are now trying to provide
more detailed conditions for the convergence rates at some extreme
cases. It is easy to see that $\gamma(\mu,\alpha,c_1)$ is
monotonically decreasing with respect to $c_1\in (0, 1]$, which
implies
$$
\frac{(1+\mu)(\mu-\alpha)}{\mu(1+\alpha)} \le \gamma(\mu,\alpha,c_1)
< \mu+1  < 2.
$$
When $c_1$ is close to $1$, i.e., $D$ is relatively small,  we have
$$
\gamma (\mu,\alpha,c_1) \approx
\frac{(\mu+1)(\mu-\alpha)}{\mu(1+\alpha)}\,.
$$
Hence for any $\theta_i$ satisfying
$$
\delta_2(1+\beta)\theta_i \le
\frac{(\mu+1)(\mu-\alpha)}{\mu(1+\alpha)} < 2
\frac{1-\alpha}{1+\alpha} < \frac{2}{\kappa_1},
$$
we know $F_i\le\mu\,I$, while for $\tilde \mu$ in the following range
$$
1-\tilde{\mu}\le \frac{\delta_1
(1-\beta)}{1+\frac{\alpha}{\beta+\alpha} \frac{2}{\kappa_1}}
\theta_i ,
$$
we know $F_i \ge - \tilde{\mu}\,I$.

In the case that $c_0$ and $c_1$ are both close to $0$,
i.e., $D$ dominates $B^tA^{-1}B$, we see
$$
\gamma(\mu,\alpha,c_1) \approx \mu+1,\quad \delta_1\approx 1\,,
~~\delta_2 \approx 1 .
$$
Thus for $\theta_i$ satisfying
\begin{equation} \label{Th2}
\begin{array}{l}
\ds \theta_i < \frac{\mu+1}{1+\beta}
\end{array} \end{equation}
then $F_i\le \mu\,I$. On the other hand, it follows from \eqref{R2} that if
\begin{equation} \label{rate2}
1-\tilde{\mu} \le \theta_i\,(1-\beta)
\end{equation}
then $F_i \ge -\tilde{\mu}\,I$.

From above we can see that
the convergence rate $\rho=\max \{\mu,\tilde{\mu} \}$ can be
estimated using Theorem~\ref{thm:rates1} and
\eqref{Th2}--\eqref{rate2} when $D$ is dominant in the approximate
Schur complement $H=B^T\hat A^{-1}B+D$.

\begin{cor}\label{cor:D0}
Let $D=0$ in \eqref{saddle}. Then for any damping parameter
$\theta_i$ satisfying
$$
\theta_i \le \frac{\lambda}{\kappa_1}\,,
$$
Algorithm \ref{algo:LinInexUzawaVariableRelax} converges; and if the
eigenvalues of $S^{1/2}Q_i^{-1}S^{1/2}$ is clustered around
$(1-\alpha)/(1+\alpha)$, the algorithm achieves approximately the
optimal rate $\sqrt{\alpha}$.
\end{cor}

\noindent {Proof}. When $D=0$,
the error propagating matrix $F_i$ in (\ref{eq:Fi}) becomes
$$
\tilde F_i=\left( \begin{array}{cc} \alpha\,(I+\Sigma_0^t U^t Q_i^{-1} U
\Sigma_0) &
\sqrt{\alpha}\, \Sigma_0^t U^tQ_i^{-1}U\Sigma_0 \\ \\
\sqrt{\alpha}\, \Sigma_0^t U^tQ_i^{-1} U \Sigma_0 & -(I-\Sigma_0^t
U^t Q_i^{-1} U \Sigma_0)\end{array} \right).
$$
Then we have
$$
\left(\begin{array}{c} E^{(1)}_{i+1} \\ E^{(2)}_{i+1}\end{array}
\right)
=\tilde{F}_i\left(\begin{array}{c}\frac{1}{\alpha}V^t(I-\omega_i
A^{\frac{1}{2}}\hat{A}^{-1}A^{\frac{1}{2}}) VE^{(1)}_i \\-E^{(2)}_i
\end{array}\right)\,.
$$
Clearly, $\tilde F_i$ is  a function of the single matrix
$R^tQ_i^{-1}R=\Sigma_0^t U^tQ_i^{-1}U\Sigma_0$. Let $M= R^t Q_i^{-1}
R$, then
$$
\tilde F_i - \mu I =\left( \begin{array}{cc} (\alpha-\mu)I+\alpha M & \sqrt{\alpha}M\\
\sqrt{\alpha}M & -(\mu+1)I+M \end{array} \right).
$$
Using the factorization
$$
 \begin{pmatrix} \sqrt{\alpha}M & -(\mu+1)I+M \\
(\alpha-\mu)I+\alpha M & \sqrt{\alpha}M  \end{pmatrix} =
 \begin{pmatrix} \sqrt{\alpha}M & 0 \\
(\alpha-\mu)I+\alpha M & X \end{pmatrix}
 \begin{pmatrix} I & Y \\ 0 & I  \end{pmatrix},
$$
where $X=\sqrt{\alpha} M - [(\alpha-\mu)I + \alpha M]
\alpha^{-\frac{1}{2}} M^{-1} [-(\mu+1)I + M]$,
$Y=\alpha^{-\frac{1}{2}} M^{-1} [-(\mu+1)I + M]$, we know that
$\det(\tilde F_i - \mu I)=0 $ is equivalent to $$\det( \alpha M^2 -
[ (\mu - \alpha) I - \alpha M] [(\mu +1) I - M] ) = 0.$$

 Let $z$ be an eigenvalue of $R^tQ_i^{-1}R$, then the corresponding
eigenvalue $\mu$ of $\tilde F_i$ satisfies
$$
f(\mu)=(\mu-\alpha(1+z))(\mu+1-z)-\alpha z^2
=\mu^2+(1-\alpha-(1+\alpha)z)\mu-\alpha=0.
$$
It is easy to see that
$
f(0)=-\alpha<0,
$
and
$
f(-1)=(1+\alpha)z>0.
$
If
$$
f(1)=2(1-\alpha)-(1+\alpha)z>0,
$$
then we know $\mu\in (-1,1)$. This is equivalent to
\begin{equation}\label{eq:z_estimate1}
z<\frac{2(1-\alpha)}{1+\alpha}=2\frac{\lambda}{\lambda_0}=\frac{2}{\kappa_1}\,.
\end{equation}
Noting that $S^{\frac 12}Q_i^{-1}S^{\frac 12}$ has the same eigenvalues
as $R^tQ_i^{-1}R$, we know from Lemma~\ref{lem:eig_estimate1} that
$$
\theta_i\frac{1-\beta}{\lambda_0} \le
z\le\theta_i\frac{1+\beta}{\lambda}\,,
$$
which indicates that condition (\ref{eq:z_estimate1}) holds if
$
\theta_i \le {\lambda}/{\kappa_1}\,.
$
This proves the first part of Corollary~\ref{cor:D0}.

To see the second part, we know if $z$ is clustered around
$$
\frac{1-\alpha}{1+\alpha}=\frac{1}{\kappa_1},
$$
$f(\mu)$ approaches $\mu^2-\alpha$, indicating that Algorithm
\ref{algo:LinInexUzawaVariableRelax} achieves approximately the
optimal convergence rate $\sqrt{\alpha}$. $\square$

\begin{rem}\label{rem:1b} A few remarks are in order.
\begin{enumerate}

\item
The convergence of Algorithm \ref{algo:LinInexUzawaVariableRelax}
was analyzed in \cite{huzou01} when $D=0$ under Assumption
\eqref{eq:huzoucond}, and the convergence was established by
evaluating the maximum eigenvalues of $\hat{F}_i^t\hat{F}_i$
directly, where $\hat{F}_i$ is a non-symmetric matrix given by
$$
\hat{F}_i=\left( \begin{array}{cc} \alpha_i\,(I+\Sigma_0^t U^t
Q_i^{-1} U \Sigma_0) &
-\sqrt{\alpha}\, \Sigma_0^t U^tQ_i^{-1}U\Sigma_0 \\ \\
\frac{\alpha_i}{\sqrt{\alpha}}\,\Sigma_0^t U^tQ_i^{-1} U \Sigma_0 &
(I-\Sigma_0^t U^t Q_i^{-1} U \Sigma_0)\end{array} \right).
$$
Note that $\hat{F}_i$ is a function of the single matrix
$\Sigma_0^t U^tQ_i^{-1}U\Sigma_0$.
Our estimate \eqref{eq:rate1}  is different from the one in
\cite{huzou01} since it contains the additional decay factor
${\alpha_i^2}/{\alpha^2}$. Moreover, a direct extension of the
analysis in \cite{huzou01} for the general case of $D\ne 0$ is
considerably difficult since the analysis in \cite{huzou01} depends
on the fact that $\hat F_i$ is  a function of a single matrix, but
the corresponding matrix ${F}_i$ for the case $D\neq0$ involves two
different matrices.

\item For the estimate of eigenvalues of
$S^{\frac{1}{2}}Q_i^{-1}S^{\frac{1}{2}}$ in
Lemma~\ref{lem:eig_estimate1} the estimate \eqref{est} may be very
conservative and can be replaced by the specific conditioning of
$\hat{A}^{-1}$ on $Range(B)$, i.e., \bb\label{eq:new_est1} \gamma_1
B^t\hat{A}^{-1}B\le B^tA^{-1}B\le \gamma_2 B^t\hat{A}^{-1}B\,. \ee
As a consequence the estimate of the range of eigenvalues of
$S^{\frac{1}{2}}Q_i^{-1}S^{\frac{1}{2}}$ is sharper and the
convergence rate can be improved. 

\item In all the above estimates $\beta$ can be replaced by
$\beta_i$ ($\beta_i\le\beta$) since eigenvalues of
$S^{\frac{1}{2}}Q_i^{-1}S^{\frac{1}{2}}$ is bounded in terms of
$\beta_i$ in Lemma~\ref{lem:eig_estimate1}.
In practice $\beta_i$ may be much smaller than
$\beta$, thus it may result in much sharper estimate for the lower
bound of the eigenvalues of $F_i$.
\end{enumerate}
\end{rem}


\section{Nonlinear Preconditioners}

Our analysis in the previous sections still applies when the
preconditioner $\hat{A}^{-1}$ for $A$ in (\ref{saddle}) is replaced
by a more general one. A general preconditioner is a nonlinear
mapping $\Psi_A: \mathbb{R}^n \to \mathbb{R}^n$ for the linear
system
$$
Ax=\xi
$$
such that $\Psi_A(\xi)$ gives an approximation of the solution $x$
with certain accuracy.
We assume that $\Psi_A$ satisfies
\begin{eqnarray}
|\Psi_A(\xi)-A^{-1}\xi|_A &\le& \delta |A^{-1}\xi|_A \quad
\forall\,\xi\in \mathbb{R}^n\,, \label{cond1}\\
|\Psi_A(Bd)-A^{-1}Bd|_A &\le& \delta_0 |A^{-1}Bd|_A\quad \forall\,
d\in \mathbb{R}^n \label{cond2}
\end{eqnarray}
for some $\delta$, $\delta_0 \in (0,1)$. General preconditioners of
this type can be realized, for example, by the approximate inverse
generated via the preconditioned conjugate gradient (PCG) iteration,
or by one sweep of a multigrid method with conjugate gradient
smoothing. With the help of this general preconditioner $\Psi_A$ we
consider the following iterative method for solving the generalized
saddle-point system (\ref{saddle}).

\medskip
\begin{algorithm}
\caption{Nonlinear inexact Uzawa algorithm when good approximate
Schur complement available. } \label{algo:NonlinInexUzawa1}
\begin{algorithmic}
\STATE{
\begin{enumerate}
\item Compute $f_i=f- A x_i - B y_i$, $r_i=\Psi_A(f_i)$, and the
relaxation parameter
$$ \omega_i =   \frac{ \langle f_i, r_i \rangle }{ \langle A r_i, r_i \rangle } ~~\text{for}~ f_i \neq 0  ~\text{($\omega_i=1$,
   otherwise)}. $$
\item Update $x_{i+1}=x_{i}+\omega_i\,\Psi_A(f-Ax_{i}-By_{i}) = x_i + \omega_i r_i $;
\item Compute $g_i=B^t x_{i+1} - D y_i -g$, $s_i=\hat{S}^{-1}g_i$, and
$$
\tau_i=\theta_i \frac{ \langle g_i,s_i \rangle }{ \langle
\Psi_A(Bs_i),Bs_i \rangle + \langle Ds_i,s_i \rangle } ~~\text{for}~
s_i \neq 0 ~(\tau_i=1, \text{otherwise});
$$
\item Update $y_{i+1}=y_{i}+\tau_i\hat{S}^{-1}(B^tx_{i+1}-Dy_{i}-g) = y_i + \tau_i s_i$.
\end{enumerate}
}
\end{algorithmic}
\end{algorithm}

Using condition (\ref{cond1}), one can find a symmetric and positive definite matrix $Q_{i,A}$
such that (cf.~\cite{s3})
$$
Q_{i,A}\Psi_A (f_i)=f_i
$$
and
$$
|I-A^{\frac 12} Q_{i,A}^{-1}A^{\frac{1}{2}}|\le \delta\,.
$$
Similarly, there exists a symmetric and positive definite matrix $Q_{i,B}$ such that
$$
Q_{i,B}\Psi_A(Bs_i)= B s_i
$$
and
$$
|B^tA^{-1}B-B^t Q_{i,B}^{-1} B|\le \delta_0.
$$
Using the same arguments as in steps (3.1)--(3.4), we can obtain
$$\begin{array}{l}
\ds e^y_{i+1}=Q_i^{-1}B^tA^{-\frac{1}{2}}(I-\omega_i A^{\frac{1}{2}}
Q_{i,A}^{-1}A^{\frac{1}{2}})A^{-\frac{1}{2}}f_i +(I-Q_i^{-1}S)e^y_i,
\\ \\
\ds A^{-\frac{1}{2}}f_{i+1}
=(I+A^{-\frac{1}{2}}BQ_i^{-1}B^tA^{-\frac{1}{2}}) (I-\omega_i
A^{\frac{1}{2}} Q_{i,A}^{-1}A^{\frac{1}{2}})A^{-\frac{1}{2}}f_i
-A^{-\frac{1}{2}}BQ_i^{-1}Se^y_i.
\end{array} $$
Let $H_i=B^t Q_{i,A}^{-1}B+D$,  and
$$
\beta_i=\frac{|(I-\hat{\tau}_iH_i^{\frac{1}{2}}\hat{S}^{-1}
H_i^{\frac{1}{2}})H_i^{-\frac{1}{2}}g_i|}{|H_i^{-\frac{1}{2}}g_i|}\,,
\quad \kappa=\mbox{cond}(\hat{S}^{-1}(B^tA^{-1}B+D)).
$$
Then one can prove that there exists $\beta=\beta(\delta_0,\kappa)$
such that $\beta_i\le \beta\le 1$ as it was done in the proof of
Lemma~\ref{lem:eig_estimate1}. Consequently, one can prove
Lemma~\ref{lem:eig_estimate1} and (\ref{eq:new_est1}) in
Remark~\ref{rem:1b} with
$\gamma_1=1-\delta_0,\;\gamma_2=1+\delta_0$, thus we can carry out
exactly the same convergence analysis as we did in the previous
sections for the nonlinear inexact Uzawa algorithm above.

When there is no good preconditioner for the Schur complement
system, especially when $\text{cond}(\hat{S}^{-1} S) \gg
\text{cond}(\hat{A}^{-1}A)$, we use a nonlinear solver, for example,
CG, to solve $H z = \zeta$, where $H=B^t \hat{A}^{-1} B + D$, and
get the approximate solution $\psi_H (\zeta)$.

\begin{algorithm}
\caption{Nonlinear inexact Uzawa algorithm when no good approximate
Schur complement available. } \label{algo:NonlinInexUzawa2}
\begin{algorithmic}
   \STATE{1. Compute $f_i = f - (A x_i + B y_i)$, $r_i = \hat{A}^{-1} f_i$, 
   and the relaxation parameter
   $$ \omega_i =   \frac{ \langle f_i, r_i \rangle }{ \langle A r_i, r_i \rangle } ~~\text{for}~ f_i \neq 0  ~\text{($\omega_i=1$,
   otherwise)}. $$
   }
   \STATE{2. Update $x_{i+1} = x_i + \omega_i r_i$.
   }
   \STATE{3. Compute $g_i = B^t x_{i+1} - D y_i -g$, $s_i = \Psi_H (g_i)$, and the parameter
   $$\tau_i =   \theta_i \frac{ \langle g_i, s_i \rangle }{ \langle Hs_i,s_i \rangle } ~~\text{for}~ s_i \neq 0 ~\text{($\tau_i=1$,
   otherwise)}. $$
   }
   \STATE{4. Update $y_{i+1} = y_i +   \tau_i s_i$.
   }
\end{algorithmic}
\end{algorithm}

Assume that $|\Psi_H (\zeta) - H^{-1} \zeta |_H \leq \delta_H |
H^{-1} \zeta |_H, \forall \zeta \in \mathbb{R}^m $. There is a
symmetric positive definite matrix $\hat{Q}_i$ (see Lemma 9 in
\cite{s3}) such that $\hat{Q}_i^{-1} g_i=\Psi_H(g_i) $ and all
eigenvalues of the matrix $\hat{Q}_i^{-1} H$ are in the interval
$[1-\delta_H, 1+\delta_H]$. That is,
$$ (1-\delta_H) \langle \hat{Q}_i \psi, \psi \rangle  \leq  \langle H \psi, \psi \rangle  \leq (1 + \delta_H)  \langle \hat{Q}_i \psi, \psi \rangle ,
\quad \forall \psi \neq 0.$$ Suppose that $  \langle  S \phi, \phi
\rangle = \lambda  \langle \hat{Q}_i \phi, \phi  \rangle $, where
$\lambda$ is the eigenvalue of $\hat{Q}_i^{-1} S$. We can verify
that
$$  \langle S \phi, \phi \rangle  =  \langle A^{-1} B \phi, B \phi \rangle  +  \langle D \phi,\phi \rangle  =  \langle \hat{A}^{\frac{1}{2}}
A^{-1} \hat{A}^{\frac{1}{2}} \hat{A}^{-\frac{1}{2}} B \phi,
\hat{A}^{-\frac{1}{2}} B \phi  \rangle  +  \langle D \phi, \phi
\rangle . $$
Let $\mu_1$ and
$\mu_2$ are the minimal and maximal eigenvalues of $\hat{A}^{-1} A$,
and we have
$$ \mu_1  \langle \hat{A}^{-\frac{1}{2}} B \phi,
\hat{A}^{-\frac{1}{2}} B \phi  \rangle  +  \langle D \phi,\phi
\rangle \leq  \langle  S \phi, \phi  \rangle  \leq  \mu_2
 \langle \hat{A}^{-\frac{1}{2}} B \phi, \hat{A}^{-\frac{1}{2}} B \phi \rangle  +  \langle D
\phi,\phi \rangle .$$
Assuming that the spectra of $\hat{A}^{-1} A$ are around 1 and
$\mu_1 \leq 1 \leq \mu_2$,
 we obtain
$$ \mu_1  \langle B^t \hat{A}^{-1} B \phi,  \phi \rangle  + \mu_1  \langle D \phi,\phi \rangle   \leq  \langle  S \phi,
\phi \rangle  \leq  \mu_2  \langle B^t \hat{A}^{-1} B \phi,  \phi
\rangle  + \mu_2  \langle D \phi,\phi \rangle .$$ That is,
$$ \mu_1  \langle H \phi, \phi \rangle  \leq \lambda  \langle  \hat{Q}_i \phi, \phi  \rangle  \leq \mu_2  \langle H \phi,
\phi \rangle .$$ 
Hence,
$$ \mu_1 (1-\delta_H)  \langle \hat{Q}_i \phi, \phi \rangle  \leq \lambda  \langle  \hat{Q}_i \phi, \phi  \rangle
\leq \mu_2 (1+\delta_H)  \langle \hat{Q}_i \phi, \phi \rangle .$$
One can directly check that $$\text{cond}(\hat{Q}_i^{-1} S ) \leq
\frac{1 +\delta_H}{1-\delta_H} \frac{\mu_2}{\mu_1} = \frac{1
+\delta_H}{1-\delta_H} \text{cond}(\hat{A}^{-1} A).$$

Therefore, the nonlinear solver $\Psi_H(g_i) (= \hat{Q}_i^{-1} g_i)$
corresponds to a new preconditioner $\hat{Q}_i$ such that
cond($\hat{Q}_i^{-1} S$) is much more improved than
cond($\hat{S}^{-1} S$) and has about the same order as
$\text{cond}(\hat{A}^{-1} A)$.
Algorithm \ref{algo:LinInexUzawaVariableRelax} can be recovered if
we replace $\hat{Q}_i$ by $\hat{S}$ in Algorithm
\ref{algo:NonlinInexUzawa2}. Obviously Algorithm
\ref{algo:NonlinInexUzawa2} can  be regarded as a variant of the
previous Algorithm \ref{algo:LinInexUzawaVariableRelax}, and similar
convergence analysis can be performed for Algorithm
\ref{algo:NonlinInexUzawa2}.

\section{Nonsymmetric case}\label{sec:nonsymm}

In this section we consider the convergence of Algorithm
\ref{algo:LinInexUzawaVariableRelax} for the case when $A$ in
(\ref{saddle}) is nonsymmetric. This study seems to be new, and
still no such investigations are available in the literature. Let
$A_0$ be the symmetric part of $A$,  with $A_0$ being positive
definite. Let
$$
J=A_0^{\frac{1}{2}}A^{-1}A_0^{\frac{1}{2}}.
$$
First, we note that the relaxation parameter $\omega_i$ in Algorithm
\ref{algo:LinInexUzawaVariableRelax} is now replaced by
$$
\omega_i=\frac{(f_i,r_i)}{(A_0r_i,r_i)}\,.
$$
Using (\ref{eq:fi}) and the iteration (\ref{eq:alg2}) for updating
$x_i$, we can write \bb \l{nons1} A_0^{\frac{1}{2}}e^x_{i+1}=
A_0^{\frac{1}{2}}(e^x_i-\omega_i\hat{A}^{-1}f_i) =(J-\omega_i
A_0^{\frac{1}{2}}\hat{A}^{-1}A_0^{\frac{1}{2}})A_0^{-\frac{1}{2}}f_i-JA_0^{-\frac{1}{2}}Be^y_i.
\ee On the other hand, using the iteration (\ref{eq:alg2}) for
updating $y_i$, the definition of $g_i$, (\ref{nons1}) and matrix
$G_i$ introduced in Lemma~\ref{lem:eig_estimate1} we derive \bb
\l{eq:long1}
\begin{array}{l}
\ds e^y_{i+1}=e^y_i-Q_i^{-1}g_i=e_i^y+Q_i^{-1}(B^te^x_{i+1}-De^y_i)
\\ \\
\ds \quad =e^y_i+Q_i^{-1}B^tA_0^{-\frac{1}{2}}((J-\omega_i
A_0^{\frac{1}{2}}\hat{A}^{-1}
A_0^{\frac{1}{2}})A_0^{-\frac{1}{2}}f_i-JA_0^{-\frac{1}{2}}Be^y_i)-Q_i^{-1}De^y_i
\\ \\
\ds\quad =Q_i^{-1}B^tA_0^{-\frac{1}{2}}(J-\omega_i
A_0^{\frac{1}{2}}\hat{A}^{-1}A_0^{\frac{1}{2}})A_0^{-\frac{1}{2}}f_i
+(I-Q_i^{-1}S)e^y_i\,,
\end{array}
\ee
where $Q_i^{-1}=\theta_iG_i^{-1}$.
Now, it follows from (\ref{eq:fi}), (\ref{nons1}) and
(\ref{eq:long1}) that
\bb \l{eq:fi+1}
\begin{array}{l}
A_0^{-\frac{1}{2}}f_{i+1}=J^{-1}A_0^{\frac{1}{2}}e^x_{i+1}+A_0^{-\frac{1}{2}}Be^y_{i+1}
\\ \\
\quad =(J^{-1}+A_0^{-\frac{1}{2}}BQ_i^{-1}B^tA_0^{-\frac{1}{2}})
(J-\omega_i
A_0^{\frac{1}{2}}\hat{A}^{-1}A_0^{\frac{1}{2}})A_0^{-\frac{1}{2}}f_i
-A_0^{-\frac{1}{2}}BQ_i^{-1}Se^y_i\,.
\end{array} \ee
Consider the singular value decomposition of matrix $B^tA_0^{-\frac{1}{2}}$,
$$
B^tA_0^{-\frac{1}{2}}=U\Sigma V^t,\quad \Sigma=[\Sigma_0, \; ~0]
$$
where $U$ is an orthogonal $m\times m$ matrix, $V$ is an
orthogonal $n\times n$ matrix, and $\Sigma_0$ is a $m\times m$
diagonal matrix with its diagonal entries being the singular values of
$B^tA_0^{-\frac{1}{2}}$.  Let $S_0=B^t A_0^{-1}B+D$ and $S_0=RR^t$.
Set
$$
E^{(1)}_i=\sqrt{\alpha}V^tA_0^{-\frac{1}{2}}f_i, \quad E^{(2)}_i=R^t
e^y_i.
$$
Using (\ref{eq:fi+1}), we have
\begin{eqnarray*}
\sqrt{\alpha} V^t A_0^{-\frac{1}{2}} f_{i+1} &=& [ V^t(J^{-1} + V X
V^t)(J-\omega_i Y)] V ( \sqrt{\alpha} V^t A_0^{-\frac{1}{2}} f_i)
\\ && - \sqrt{\alpha} V^t V \Sigma U^t Q_i^{-1} S S_0^{-1} R (R^t
e_i^y),
\end{eqnarray*}
where $X=\Sigma^t U^tQ_i^{-1}U\Sigma$ and
$Y=A_0^{\frac{1}{2}}\hat{A}^{-1}A_0^{\frac{1}{2}}$. Noticing that
$(V^t J^{-1} + X V^t) (J - \omega_i Y)= (I+X) V^t (I-\omega_i Y) + X
V^t(J-I) - \omega_i V^t(J^{-1}-I)Y$ and $S S_0^{-1} R= (S-S_0 +
S_0)S_0^{-1} R = R - (S_0-S)S_0^{-1} R = R-(S_0-S)R^{-t}$, we
rewrite the formula above as follows,
\begin{eqnarray*}
E^{(1)}_{i+1} &=& [(I+X) V^t (I-\omega_i Y) + X V^t(J-I) - \omega_i
V^t(J^{-1}-I)Y] \alpha (I-\omega_i Y)^{-1} V  \\
& &  \cdot \frac{1 } {\alpha}  V^t (I-\omega_i Y) V E^{(1)}_i
  - \sqrt{\alpha} V^t V \Sigma U^t Q_i^{-1}[R-(S_0-S)R^{-t}] E^{(2)}_i \\
&=& \alpha \{(I+X)   + [X V^t(J-I) - \omega_i
V^t(J^{-1}-I)Y]  (I-\omega_i Y)^{-1} V \}  \\
& &  \cdot \frac{1 } {\alpha}  V^t (I-\omega_i Y) V E^{(1)}_i
  - \sqrt{\alpha} \Sigma U^t Q_i^{-1}[R-(S_0-S)R^{-t}] E^{(2)}_i .
\end{eqnarray*}
Using (\ref{eq:long1}), we obtain
\begin{eqnarray*}
R^t e_{i+1}^y &=& \sqrt{\alpha} R^t Q_i^{-1} U \Sigma V^t
(J-\omega_i Y) \frac{1}{\alpha} V ( \sqrt{\alpha} V^t
A_0^{-\frac{1}{2}} f_i ) \\
&&  + R^t (I- Q_i^{-1} S ) R^{-t} ( R^t e_i^y ).
\end{eqnarray*}
Noticing that $(J-\omega_i Y) V = (I - \omega_i Y + J -I)(I -
\omega_i Y)^{-1} V V^t (I - \omega_i Y) V = [V + (J-I)(I - \omega_i
Y)^{-1} V] V^t (I - \omega_i Y) V$, and $ R^t(I-Q_i^{-1} S) R^{-t} =
R^t [I - Q_i^{-1}(S_0 + S - S_0)] R^{-t} = I - R^t Q_i^{-1} R - R^t
Q_i^{-1} (S-S_0) R^{-t} $, we rewrite the formula above as follows,
\begin{eqnarray*}
E^{(2)}_{i+1} &=& \sqrt{\alpha} [R^t Q_i^{-1} U \Sigma  + R^t
Q_i^{-1} U \Sigma V^t (J-I)(I-\omega_iY)^{-1} V]
\frac{1}{\alpha} V^t (I-\omega_i Y) V E^{(1)}_i  \\
&&  - [-(I - R^t Q_i^{-1} R ) +  R^t Q_i^{-1} (S-S_0) R^{-t} ]
E^{(2)}_i.
\end{eqnarray*}


The error propagation can be reformulated as
$$\begin{array}{l}
\ds \left(\begin{array}{c} E^{(1)}_{i+1} \\ \\ E^{(2)}_{i+1}
\end{array} \right)= \left[\left(\begin{array}{cc} \alpha(I+\Sigma^t
U^tQ_i^{-1}U\Sigma) & \sqrt{\alpha}
\Sigma^t U^tQ_i^{-1}R  \\ \\
\sqrt{\alpha}R^tQ_i^{-1}U\Sigma &
-(I-R^tQ_i^{-1}R)
\end{array} \right) +\Delta\right]
\\ \\
\qquad\qquad\qquad \cdot \left( \begin{array}{c} \frac{1}{\alpha}V^t
(I-\omega_i A_0^{\frac{1}{2}}\hat{A}^{-1}A_0^{\frac{1}{2}})V E^{(1)}_i \\
\\ -E^{(2)}_i
\end{array}\right)
\end{array} $$
where $\Delta = \begin{pmatrix} \Delta_{11} & \Delta_{12} \\
\Delta_{21} & \Delta_{22} \end{pmatrix}$ with the blocks defined by
$$\begin{array}{l}
\ds \Delta_{11}=\alpha\, [ \Sigma^t U^t Q_i^{-1}U\Sigma V^t (J-I)
-\omega_i V^t
(J^{-1}-I)A_0^{\frac{1}{2}}\hat{A}^{-1}A_0^{\frac{1}{2}} ]
(I-\omega_i A_0^{\frac{1}{2}}\hat{A}^{-1}A_0^{\frac{1}{2}})^{-1}V\,,
\\ \\
\ds \Delta_{12}=-\sqrt{\alpha} \Sigma^t U^tQ_i^{-1}(S_0-S) R^{-t}\,,
\\ \\
\ds \Delta_{21}=
\sqrt{\alpha}\,R^t Q_i^{-1}U\Sigma V^t (J-I)
(I-\omega_i A_0^{\frac{1}{2}}\hat{A}^{-1}A_0^{\frac{1}{2}})^{-1}V\,,
\\ \\
\ds \Delta_{22}=R^tQ_i^{-1}(S-S_0)R^{-t}\,.
\end{array} $$
Now it follows from Theorem \ref{thm:spectral_f1} that Algorithm
\ref{algo:LinInexUzawaVariableRelax} will converge when
$$
|J-I|\,, \quad |J^{-1}-I|,\; \quad |S-S_0|
$$
are sufficiently small.

\section{Applications}\label{sec:appls}
The saddle-point system \eqref{saddle} arises from many
applications. We present a few such examples in this section.

The first example arises naturally from the standard quadratic
constrained programming with linear constraints:
\begin{equation} \label{const}
\min_{x\in \mathbb{R}^n} ~J(x)=\frac{1}{2}(Ax,x)-(f,x)
\quad\mbox{subject to} \quad Bx=g\,.
\end{equation}
If we apply the Lagragian multiplier approach with penalty for the minimization
problem (\ref{const}),
we come to solve system \eqref{saddle} for the primal variable $x$
and the Lagrange multiplier $y$,
with $D=\epsilon\,\hat{D}$, where $\epsilon>0$ is usually
a small parameter and $\hat D$ is an appropriately selected
symmetric and positive definite matrix. If we apply the above approach
iteratively with respect to $\epsilon$, then the parameter
$\epsilon$ needs not be too small.

The second example is related to the mixed formulation
for the second order elliptic equation,
$-\nabla\cdot(\mu \nabla u)+cu=f$. In some applications the flux
$p=\mu\nabla u$ is an important quality to know. For the purpose,
we may introduce the new variable $p=\mu\nabla u$,
then the elliptic equation can be written as the system
$$\begin{array}{lll}
\ds \frac{1}{\mu}\,p-\nabla u=0\,, \quad
\ds -\nabla\cdot p+cu=f\,.
\end{array} $$
When we apply the mixed finite element formulation to the above
system, we obtain a discrete system of form \eqref{saddle}.

The third example comes from the linear elasticity equation
\begin{equation} \label{elas}
-\mu\,\Delta u-\nabla((\lambda+\mu)\nabla\cdot u)=f
\end{equation}
where $\mu,\;\lambda$ are Lame coefficients. If one needs to
follow the compressiveness  of the displacement more closely,
one may introduce
a new variable $p=(\lambda+\mu)\nabla \cdot u$,
then \eqref{elas} can be equivalently written as
\begin{equation} \label{elas2}
\ds -\mu \,\Delta u-\nabla p = f\,, \quad \ds \nabla \cdot
u-\frac{1}{\lambda+\mu}p=0\,.
\end{equation}
This formulation allows us to develop some stable numerical methods
for the nearly incompressible case, $\lambda \gg 1$. Now the
application of the mixed finite element formulation to the above
system results in a discrete system of form \eqref{saddle}.

The next example arises from the following elliptic interface
problem
$$\begin{array}{l}
\ds- \nabla \cdot(\mu \nabla u)=f \quad \mbox{in} \quad \Omega\,;
\\ \\
\ds [\mu\frac{\partial u}{\partial\nu}] +\alpha\, u=g \quad \mbox{on}
\quad \Gamma\,,
\end{array} $$
where $\Omega$ is occupied by, e.g.,  two different fluids or
materials $\Omega_1$ and $\Omega_2$,  with different physical
property $\mu$ and a common interface $\Gamma=\bar\Omega_1\cap
\bar\Omega_2$. $[\mu{\partial u}/{\partial\nu}]$ stands for the jump
of the flux $\mu{\partial u}/{\partial\nu}$ across the interface. In
some applications, the jump of the flux $[\mu{\partial
u}/{\partial\nu}]$ can be an important physical quantity to know.
For this purpose, we may introduce a new variable
$p=-[\mu\frac{\partial}{\partial\nu}u]$, then the above interface
system can be written as
\begin{eqnarray*}
 \ds - \nabla \cdot(\mu \nabla u)+\gamma^* p &=& f \quad \mbox{in} \quad
\Omega\,;\\ \\
\ds \gamma u-\frac{1}{\alpha}p &=& g  \quad \mbox{on} \quad
\Gamma\,,
\end{eqnarray*}
where $\gamma$ is the trace operator from $H^1(\Omega)$ to $L^2(\Gamma)$
and $\gamma^*p \in H^1(\Omega)^*$
is defined by
$$
\langle \gamma^* p,\phi\rangle=(p,\gamma \phi)_{L^2(\Gamma)}, \quad
\forall\,p\in L^2(\Gamma), \phi\in H^1(\Omega)\,.
$$
The advantage of this formulation is that it can be easily utilized
in the domain decomposition approach for a wide class of interface
problems, e.g., one uses a subdomain solver, given the boundary
value $g$ and solves the Schur complement system
(Neumann-to-Dirichlet map) that equates the continuity of the
solution at $\Gamma$.

In addition, \eqref{saddle} can be regarded as a regularization of
the simplified saddle-point problem where the (2,2) diagonal block
vanishes, with ${D}$ arising from the regularization on $y$. This
regularization is often used to remedy the lack of the inf-sup
condition and prevent the locking phenomena; see \cite{s6, s14,
GloLeta89},
for example, the stabilized Q1-P0 finite element method on the
steady-state Stokes problem: \begin{equation}  \label{Stokes_PDE} -
\nu \Delta u + \nabla p = 0 \,, \quad \ds - \nabla \cdot u   =0
\quad \mbox{in} \quad \Omega \end{equation}
with Dirichlet boundary conditions on $\partial \Omega$, where $u$
stands for the velocity field and $p$ denotes the pressure.

\section{Numerical experiments}

In the following we present some numerical experiments to show the
performance of Algorithm \ref{algo:LinInexUzawaVariableRelax} with
parameters $\omega_i$ and $\tau_i$ selected by (\ref{eq:omega1}) and
(\ref{tau2}). As our first testing example, we consider the
two-dimensional elasticity problem \eqref{elas} and its mixed
formulation (\ref{elas2}) in the domain $\Omega=(0,1)\times (0,1)$.
For convenience we use $(u, v)$ and $(f, g)$ below to stand
respectively for the displacement vector $u$ and forcing vector $-f$
in (\ref{elas2}). The system \eqref{elas} is complemented by the
following boundary conditions
\begin{eqnarray}
&&u=0,\quad v_x=0 \quad\mbox{on}\quad x=0,\;1 .\\
&&u_y=0,\quad v=0 \quad \mbox{on}\quad y=0,\;1 .
\end{eqnarray}
We partition the domain $\Omega$ into $n^2$ equal rectangular
elements, and the displacement components $u$ and $v$ and the
pressure $p$ are approximated respectively at the staggered grids as
follows:
\begin{eqnarray}
\ds p^{i,j} &\approx & p((i-\frac{1}{2})\,h,(j-\frac{1}{2})\,h)
\quad \mbox{for $1\le i \le n$,
$1\le j\le n$}\,,\\
\ds u^{i,j-\frac{1}{2}} &\approx& u(i\,h,(j-\frac{1}{2})\,h)
\quad \mbox{for $0 \le i \le n$, $1 \le j\le n$}\,,\\
v^{i-\frac{1}{2},j} &\approx& v((i-\frac{1}{2})\,h,j\,h) \quad
\mbox{for $1 \le i \le n$, $0\le j\le n$,}
\end{eqnarray}
with the meshsize $\ds h={1}/{n}$. Applying the central difference
approximation to (\ref{elas2}) results in the following scheme:
\begin{eqnarray*}
&&\mu \frac{u^{i+1,j-\frac{1}{2}}-2u^{i,j-\frac{1}{2}}+
u^{i-1,j-\frac{1}{2}}}{h^2}+\mu
\frac{u^{i,j+\frac{1}{2}}-2u^{i,j-\frac{1}{2}}+
u^{i,j-\frac{3}{2}}}{h^2}\\
&+&\frac{p^{i+\frac{1}{2},j-\frac{1}{2}}
-p^{i-\frac{1}{2},j-\frac{1}{2}}}{h}=f^{i,j-\frac{1}{2}}\,,
\\
&&\mu \frac{v^{i+\frac{1}{2},j}-2v^{i-\frac{1}{2},j}+
v^{i-\frac{3}{2},j}}{h^2}+\mu
\frac{v^{i-\frac{1}{2},j+1}-2v^{i-\frac{1}{2},j}+
v^{i-\frac{1}{2},j-1}}{h^2}\\
&+&\frac{p^{i-\frac{1}{2},j+\frac{1}{2}}
-p^{i-\frac{1}{2},j-\frac{1}{2}}}{h}=g^{i-\frac{1}{2},j}\,,
\\
&&\frac{u^{i,j-\frac{1}{2}}-u^{i-1,j-\frac{1}{2}}}{h}+
\frac{v^{i-\frac{1}{2},j}-u^{i-\frac{1}{2},j-1}}{h}
-\frac{1}{\mu+\lambda^{i-\frac{1}{2},j-\frac{1}{2}}}\,p^{i-\frac{1}{2},j-\frac{1}{2}}=0.
\end{eqnarray*}
Equivalently the matrices $A$, $B$ and $D$ in (\ref{saddle}) can be written
as
$$
A=\left(\begin{array}{cc} A_1 & 0\\ \\ 0& A_2\end{array}\right), \quad
B=\left(\begin{array}{c} B_1\\ \\ B_2\end{array}\right), \quad
D=diag(\frac{1}{\mu+\lambda^{i-\frac{1}{2},j-\frac{1}{2}}})\,,
$$
where
$$
A_1=I \otimes H_1+H_2 \otimes I\,,\quad A_2=I \otimes H_2+ H_1
\otimes I\,,
$$
$$
B_1=I \otimes D,\quad B_2= D \otimes I,
$$
and the tridiagonal matrices $H_1\in \mathbb{R}^{(n-1)\times (n-1)}$
and $H_2\in \mathbb{R}^{n\times n}$ are given by
$$
H_1=\left(\begin{array}{ccccc} 2 &-1 & & &
\\ -1 & 2& -1 & &  \\
& \ddots & \ddots &\ddots & \\
 & & -1 & 2 & -1 \\
& & &-1 & 2 \end{array}\right)\,, \quad
H_2=\left(\begin{array}{ccccc} 1 &-1 & & &
\\ -1 & 2& -1 & & \\
&  \ddots &  \ddots & \ddots  & \\
&  & -1 & 2 & -1 \\
& & &-1 & 2 \end{array}\right)\,.
$$
We will choose the following set of parameters in our test:
$f=0$, $g=1$ and $\mu=1$. The parameter $\lambda$ is taken to be
discontinuous:
$\lambda=1000$ in $(0.25, 0.75)\times (0.25, 0.75)$, and
$\lambda=0$ otherwise.

We have tested Algorithm \ref{algo:LinInexUzawaVariableRelax}, with
preconditioner $\hat{A}$ taken to be the Jacobi preconditioner (a
simple but poor preconditioner) and the incomplete Cholesky
factorization (\textsc{Matlab} function {\sf cholinc}   with drop
tolerance of $10^{-3}$ and no fill-in). For the Schur complement
$S$, we take the diagonal preconditioner  $\hat S=I+D$ (a simple but
poor preconditioner). Table~\ref{elstnum1}   summarizes the
convergence of Algorithm \ref{algo:LinInexUzawaVariableRelax} for
this symmetric case, where  `Iter' stands for the iteration numbers.
The first 4 columns are for the poor Jacobi preconditioner and show
numbers of iterates and CPU time (seconds) to achieve the error
$|(f_i, g_i)|<10^{-4}$ for $n=20$. The next 4 columns are for the
more reasonable preconditioner generated by the incomplete Cholesky
factorization with no fill-in for the cases $n=20$ and $n=50$. The
next 4 columns are for the good preconditioner by incomplete
Cholesky factorization with drop tolerance of $10^{-3}$ for the
cases $n=20,50,100,200$, with a total number of degrees of freedom
being $120,000$ for $n=200$. The last column is for the case when
the exact preconditioner for $A$ is used. From our experiments and
observations, the number of iterations is insensitive to mesh
refinements if good preconditioners are used. With the poor Jacobi
preconditioner, Algorithm \ref{algo:LinInexUzawaVariableRelax}
always converges. We have tested the algorithm with the damping
factor $\theta$ selected from the range $[0.01, \; 1.0]$, and
observed the convergence of the algorithm for all the cases. But for
the well-conditioned case for $A$, $\theta_i=1$ produces the best
results.  For the very ill-conditioned preconditioner $\hat{A}$,
$\theta_i$ may need to be small.

\newsavebox{\mytablebox}    
\begin{lrbox}{\mytablebox}

\begin{tabular}{c rrrr crrrr crrrr cr}
\hline  &  \multicolumn{4}{c}{Jacobi} & & \multicolumn{4}{c}{no
fill-in
Cholesky} & & \multicolumn{4}{c}{cholinc($A, 10^{-3}$) } & & Exact \\
\cline{2-5} \cline{7-10} \cline{12-15} \cline{17-17}
 $\theta_i$  & .03 & .1 & .5 & 1 & & 1  & 1 & .1 & .05 & & 1 & 1 & 1 & .1 & & 1 \\

Iter  & 659 & 737 & 906 & 1074  & & 95 & 752 & 463 & 434 & & 11 & 17  & 61 & 152 & &  5 \\

CPU   & .52 & .56 & .66 & .74 & & .22 & 9.57 & 5.70 & 5.37 & & .04 & 1.86 & 4.91 & 56.9 & & 5.62 \\

$n$  & 20 & 20 & 20 & 20 & & 20 & 50 & 50 & 50 & & 20 & 50 & 100 & 200 & & 200 \\
\hline
\end{tabular}

\end{lrbox}

 \begin{table}
 \centering \scalebox{0.9}{\usebox{\mytablebox}}    
 \caption{\label{elstnum1} Number of iterates with different $\theta_i$'s and preconditioners for linear elasticity problem. }
 \end{table}


Next we consider the Stokes flow in a rectangular domain $\Omega =
(0,1) \times (0,1)$.
 Here Dirichlet boundary conditions are used: $u =1$, $v =0$ on the top ($y=1$); $u =v =0$
on the other three sides (i.e., $x=0, x=1$, and $y=0$).
We discrete the computation domain with $Q_1 -P_0$ element, where
the velocity is located on the node, the pressure is constant in the
center of each element, and the cell width is $h=1/n$. After
discretization of (\ref{Stokes_PDE}), we obtain
\begin{equation}
\begin{bmatrix} A_0 & 0& B_1^T \\ 0 & A_0 & B_2^T \\ B_1 & B_2 & -D
\end{bmatrix}
\begin{bmatrix} u  \\ v \\ p \end{bmatrix}
= \begin{bmatrix} f_1 \\ 0 \\ 0 \end{bmatrix}  , \label{Stokes_Q1P0}
\end{equation}
where $u $, $v$ and $p$ are numbered from left to right and from
bottom to top. The coefficient matrix can be given in detail as
follows,
\begin{equation*}
\begin{bmatrix} {\nu}/{6} \left(M\otimes K + K \otimes M \right) & 0 & {h}/{2} \left( H_n^T \otimes H_o^T \right)
\\ 0 & {\nu}/{6} \left(M\otimes K + K \otimes M \right) & {h}/{2}
\left( H_o^T \otimes H_n^T \right)
\\ {h}/{2} \left( H_n \otimes H_o \right) & {h}/{2}
\left( H_o \otimes H_n \right) & -\beta h^2 (I \otimes T_N + T_N
\otimes I )
\end{bmatrix}  .
\end{equation*}
Here we define $ A_0 = {\nu}/{6} \left(M\otimes K + K \otimes M
\right)$, $B_1 = {h}/{2} \left( H_n \otimes H_o \right)$, $ B_2=
{h}/{2} \left( H_o \otimes H_n \right)$, $D = \beta h^2 (I \otimes
T_N + T_N \otimes I )$,
where $ M=\mathrm{tridiag}(1,4,1)\in \mathbb{R}^{(n-1) \times
(n-1)}$,
$K= \mathrm{tridiag}(-1,2,-1)\in \mathbb{R}^{(n-1) \times (n-1)}$,
$T_N= \mathrm{tridiag}(-1,2,-1) - e_1 e_1^T -e_n e_n^T \in
\mathbb{R}^{n \times n}$,
and $H_o$, $H_n$ are bidiagonal matrices with $H_o = {\sf
sparse}(1:n-1,1:n-1,-{\sf ones}(1,n-1),n,n-1) + {\sf sparse}(2:n
,1:n-1,  {\sf ones}(1,n-1),n,n-1) \in \mathbb{R}^{n \times (n-1)}$,
$H_n = {\sf sparse}(1:n-1,1:n-1, {\sf ones}(1,n-1),n,n-1) + {\sf
sparse}(2:n  ,1:n-1, {\sf ones}(1,n-1),n,n-1) \in \mathbb{R}^{n
\times (n-1)}$. Here {\sf sparse} and {\sf ones} are \textsc{Matlab}
notations, $e_1$ and $e_n$ are the first and $n$-th column vector of
unit matrix $I_n$.
For the right hand side, $ f_1 = \left( 6 \times
\frac{\nu}{6}\right) \left( \epsilon_{n-1} \otimes \epsilon \right)
\in \mathbb{R}^{(n-1)^2
  \times 1} $, where  $\epsilon_{n-1}$ is the ($n-1$)th
column vector of unit matrix $I_{n-1}$, and $\epsilon =[1, \cdots,
1]^T \in \mathbb{R}^{(n-1) \times 1}$.  The choice of $\beta$
represents a trade-off between stability and accuracy.  We use
$\beta=0.25$ for the local stabilization and $\beta=1$ for the
global stabilization.
%
The iteration stops when the residual 
$\max\{||f_i||, ||g_i||\} < 10^{-6}$. The iteration numbers and
computation times are listed in Table \ref{tab:Stokes}. We compare
the iteration numbers for using different preconditioners. The
preconditioner for $A$ include Jacobi
 iteration, the incomplete Cholesky decomposition with no fill-in or with
tolerance $10^{-3}$, and the exact solver as well. The
preconditioner for Schur complement is the pressure mass matrix for
all cases. The CPU times (in seconds) are given correspondingly.

\begin{table}[!htbp]
\begin{center}
\begin{tabular}{c c l cccccccc}\hline
 && & \multicolumn{2}{c}{$\theta=0.5$} & \multicolumn{2}{c}{$\theta=0.3$} & \multicolumn{2}{c}{$\theta=0.1$} &  \multicolumn{2}{c}{$\theta=0.05$}\\ \cline{4-11}
 && & Iter & CPU & Iter & CPU & Iter & CPU & Iter & CPU \\ \hline
 \multirow{8}{*}{ \begin{sideways} $\nu=1$ \end{sideways} }&\multirow{4}{*}{ \begin{sideways} $n=32$ \end{sideways} }
 & Jacobi           & 2006 & 0.57 & 891 & 0.26  & 725 & 0.21 & 749 & 0.21 \\
  && cholinc(`0')   & 192 & 0.081 & 164 & 0.069 & 139 & 0.061 & 156 & 0.065 \\
   && cholinc($10^{-3}$) &37 & 0.022 & 47  & 0.028 & 93 & 0.056 & 175 & 0.11 \\
    && Exact            & 37 & 0.27 & 45 & 0.32 & 98 & 0.71 & 184 & 1.31 \\
    \cline{2-11}
 &\multirow{4}{*}{ \begin{sideways} $n=64$ \end{sideways} }
 & Jacobi & 16823 & 17.4 & 14518 & 15.1  & 3329 & 3.50 & 2845 & 3.02 \\
  && cholinc(`0') & 873 & 1.51 & 779 & 1.37 & 494 & 0.87 & 343 & 0.62 \\
   && cholinc($10^{-3}$) & 38 & 0.12 & 55  & 0.17 & 80 & 0.25 & 147 & 0.46 \\
    && Exact            & 36 & 1.70 & 48 & 2.25 & 94 & 4.52 & 177 & 8.34 \\
    \hline
 \multirow{8}{*}{ \begin{sideways} $\nu=0.01$ \end{sideways} }&\multirow{4}{*}{ \begin{sideways} $n=32$ \end{sideways} }
 & Jacobi       & 4103 & 1.19 & 1318 & 0.38  & 1278 & 0.37 & 1300 & 0.38 \\
  && cholinc(`0') & 295 & 0.13 & 203 & 0.094 & 235 & 0.097 & 291 & 0.12 \\
   && cholinc($10^{-3}$) & 101 & 0.061 & 117  & 0.071 & 169 & 0.10 & 271 & 0.16 \\
    && Exact & 80 & 0.57 & 115 & 0.85 & 169 & 1.21 & 269 & 1.96 \\
    \cline{2-11}
 &\multirow{4}{*}{ \begin{sideways} $n=64$ \end{sideways} }
 & Jacobi & 22026 & 23.1 & 3884 & 4.06  & 2777 & 2.91 & 3756 & 3.92 \\
  && cholinc(`0') & 1385 & 2.37 & 755 & 1.30 & 391 & 0.67 & 386 & 0.67 \\
   && cholinc($10^{-3}$) & 143 & 0.45 & 117 & 0.37 & 160 & 0.50 & 242 & 0.75\\
    && Exact & 77 & 3.60 & 95 & 4.47 & 151 & 7.14 & 247 & 11.5 \\
    \hline
\end{tabular}
\caption{\label{tab:Stokes} Stokes problem. }
\end{center}
\end{table}



%
The third testing case is a purely algebraic example from
\cite{LuZhang_SIMAX10}. Consider the linear system \eqref{saddle}
with $A=(a_{ij})_{n \times n}$, $B = [T; 0] \in \mathbb{R}^{n \times
m}$, and $D=I$, where
$$ a_{ij} = \frac{1}{\sqrt{2 \pi} \sigma} e^{\frac{-|i-j|^2}{2\sigma^2}}, \quad T = \frac{1}{1000}
\text{tridiag} (1,4,1) \in \mathbb{R}^{m \times m}.$$ We set $\sigma
= 1.5$. The right hand side is chosen such that the exact solution
is a vector of all ones. Note that $A$ is an ill-conditioned
Toeplitz matrix. Fortunately, the Schur complement $S$ is
well-conditioned for $n=800$ and $m=600$, or $n=1600$ and $m=1200$. 
 We set $\hat{S} = 2I$ as the preconditioner.
%



\begin{table}[!htbp]
\begin{center}
\begin{tabular}{lccccccccccc}\hline
 & \multicolumn{5}{c}{$n=800, m=600$} & & \multicolumn{5}{c}{$n=1600, m=1200$} \\
 \cline{2-6} \cline{8-12}
\multirow{2}{*}{$\theta_i$} & \multicolumn{2}{c}{Iter} & &
\multicolumn{2}{c}{CPU} & & \multicolumn{2}{c}{Iter} & &
\multicolumn{2}{c}{CPU} \\  \cline{2-3} \cline{5-6} \cline{8-9}
\cline{11-12}
 & Jacobi & Exact & & Jacobi & Exact & &  Jacobi & Exact & &  Jacobi &
 Exact\\
\hline

0.05& 263 & 263  & & 1.10 & 30.2 & & 263 & 263 & & 3.69 & 129.0 \\
0.1 & 206 & 129  & & 0.87 & 14.9  & & 129 &129 & & 1.86 & 63.1  \\
0.5 & 171 & 21  & & 0.72 & 2.53  & & 150 & 21 & & 2.14 & 10.3  \\
0.9 & 183 & 7 &  & 0.82 & 0.83    & & 143 & 7  & & 2.07 & 3.44 \\

\hline
\end{tabular}
\caption{\label{tab:AlgebraicCoef} The purely algebraic example.
}
\end{center}
\end{table}

As our last testing example, we consider the nonsymmetric
saddle-point system (\ref{saddle}) arising from the discretization
of the mixed formulation of the following system
$$
-\mu\,\Delta u+b\left(\begin{array}{c} \frac{\partial u_1}{\partial x_1} \\  \\
\frac{\partial u_2}{\partial x_2}\end{array}\right)+\nabla p=f\,,
$$
which is a compressible linearized Navier-Stokes system.
Numerical results are summarized in Table~\ref{elstnum2}.

\begin{table}

\begin{tabular}{c rrrrr c rrr c rrr}  \hline
&  \multicolumn{5}{c}{no fill-in Cholesky} & &  \multicolumn{3}{c}{cholinc($A,10^{-3}$)} & & \multicolumn{3}{c}{Exact} \\
\cline{2-6} \cline{8-10} \cline{12-14}
$\theta_i$ & .05 &  .05 &  .05 &  .05 &  .05  & & .03 & 1 & 1 & & .03 & 1 & 1 \\

Iter & 343 & 315 & 355 & 438 & 431 & & 1122 & 33 & 30 & & 660 & 21 & 20  \\

CPU& 4.23 & 3.92 & 4.59 & 5.45 &  5.35  & & 18.1 & .65 & .54 & & 796.2 & 21.2 & 20.5   \\

$b$ & 40 & 20 & 10 & 4 & 2 & & 10 & 4 & 2  & &  10 & 4 & 2   \\
\hline
\end{tabular}
\caption{\label{elstnum2} Nonsymmetric case with $n=50$ and
different $b$'s. }
\end{table}

The first five columns are for the preconditioner generated by the
incomplete Cholesky factorization with no fill-in for the case
$n=50$. The next three columns are for the preconditioner by the
incomplete Cholesky factorization with drop tolerance of $10^{-3}$
for $n=50$. The last three columns are for the case with exact
preconditioner for $A$ with $n=50$. The number of iterations depends
significantly on $b$ (the magnitude of the convection term). The
algorithm may fail to converge when $|b|$ is very large, which is
consistent with the convergence analysis in
Section~\ref{sec:nonsymm} as the symmetric part of block $A$ is not
dominant.

\end{document}